\documentclass[12pt]{article}
\usepackage{comment}
\usepackage{latexsym}
\usepackage{amssymb}
\usepackage{graphicx}
\usepackage{color}
\usepackage{amsmath}
\usepackage{ngerman}
\usepackage{float}
\usepackage{mathrsfs} %fonts veramente calligrafici con \mathscr{B}

\newtheorem{theorem}{Theorem}

\newtheorem{corollary}[theorem]{Corollary}

\newtheorem{definition}[theorem]{Definition}

\newtheorem{lemma}[theorem]{Lemma}

\newenvironment{proof}[1][Proof]{\textbf{#1.} }{\ \rule{0.5em}{0.5em}}
\newcommand{\kom}[1]{}
\renewcommand{\kom}[1]{{\bf [#1]}}
\definecolor{blau}{rgb}{0.1,0.0,0.9}

\newcounter{komcounter}
\numberwithin{komcounter}{section}

\begin{document}
\selectlanguage{english}
%\title{Estimates of $p$-harmonic measures and Phragmen-Lindel\"of theorems in planar cones}

%\title{Phragmen-Lindelöf theorems and p-harmonic measure in planar sectors}

\title{Estimates of $p$-harmonic functions in planar sectors}

\author{
%NIKLAS L. P. LUNDSTR\"{O}M{\small{$^1$}},
%JESPER SINGH{\small{$^{1}$}}
Niklas L.P. Lundström,
Jesper Singh
\\\\
%{\small{$^1$}}
\it \small Department of Mathematics and Mathematical Statistics, Ume{\aa} University,\\
\it \small SE-90187 Ume{\aa}, Sweden\/{\rm ;}
\it \small jesper.singh@umu.se; niklas.lundstrom@umu.se
}

\maketitle

%%%%%%%%%%%%%%%%%%%%%%%%%%%%%%%%%%%%%%%%%%%%
%%%%%%%%%%%%%%%%%%%%%%%%%%%%%%%%%%%%%%%%%%%%
%%%%%%%%%%%%%%%%%%%%%%%%%%%%%%%%%%%%%%%%%%%%
%%%%%%%%%%%%%%%%%%%%%%%%%%%%%%%%%%%%%%%%%%%%
%%%%%%%%%%%%%%%%%%%%%%%%%%%%%%%%%%%%%%%%%%%%
%%%%%%%%%%%%%%%%%%%%%%%%%%%%%%%%%%%%%%%%%%%%
%%%%%%%%%%%%%%%%%%%%%%%%%%%%%%%%%%%%%%%%%%%%
%%%%%%%%%%%%%%%%%%%%%%%%%%%%%%%%%%%%%%%%%%%%
%%%%%%%%%%%%%%%%%%%%%%%%%%%%%%%%%%%%%%%%%%%%
%%%%%%%%%%%%%%%%%%%%%%%%%%%%%%%%%%%%%%%%%%%%

\begin{abstract}
\noindent
Suppose that $p \in (1,\infty]$, $\nu \in [1/2,\infty)$,
$\mathcal{S}_\nu = \left\{ (x_1,x_2) \in  \mathbb{R}^2 \setminus \{(0, 0)\}: |\phi| <
\frac{\pi}{2\nu}\right\}$, where $\phi$ is the polar angle of $(x_1,x_2)$.
Let $R>0$ and $\omega_p(x)$ be the $p$-harmonic measure of $\partial B(0,R) \cap \mathcal{S}_\nu$ at $x$ with respect to $B(0, R)\cap \mathcal{S}_\nu$.
We prove that there exists a constant $C$ such that
\begin{align*}
C^{-1}\left(\frac{|x|}{R}\right)^{k(\nu,p)} \,  \leq \omega_p(x) \, \leq C \left(\frac{|x|}{R}\right)^{k(\nu,p)}
\end{align*}
whenever $x\in B(0,R) \cap \mathcal{S}_{2\nu}$ and
where the exponent $k(\nu,p)$ is given explicitly as a function of $\nu$ and $p$.
Using this estimate we derive local growth estimates for $p$-sub- and $p$-superharmonic functions in planar domains which are locally approximable by sectors, e.g., we conclude bounds of the rate of convergence near the boundary where the domain has an inwardly or outwardly pointed cusp.
Using the estimates of $p$-harmonic measure we also derive a sharp Phragmen-Lindel\"of theorem for
$p$-subharmonic functions in the unbounded sector $\mathcal{S}_\nu$.
Moreover, if $p = \infty$ then the above mentioned estimates extend from the setting of two-dimensional sectors to cones in $\mathbb{R}^n$.
Finally,  when $\nu \in (1/2, \infty)$ and $p\in (1,\infty)$ we prove uniqueness  (modulo normalization) of positive  $p$-harmonic functions in $\mathcal{S}_\nu$ vanishing on $\partial\mathcal{S}_\nu$.
\\


% 49L25: Viscosity solutions,
%
% 35B50: Maximum principles
%
%• 35B40: "Asymptotic behavior of solutions to PDEs" (MSC2020)
%
%• 35J92: "Quasilinear elliptic equations with p-Laplacian" (MSC2020)
%
%• 35B53: "Liouville theorems and Phragmén-Lindelöf theorems in context of PDEs" (MSC2020)
%
%• 35J25: "Boundary value problems for second-order elliptic equations" (MSC2020)
%
%• 35J60: "Nonlinear elliptic equations" (MSC2020)
%
% Secondary
%
%• 35J70: "Degenerate elliptic equations" (MSC2020)
%
%• 35D40: "Viscosity solutions to PDEs" (MSC2020)
%
%• 35B50: "Maximum principles in context of PDEs" (MSC2020)
%
%• 35J70: "Degenerate elliptic equations" (MSC2020)
%
%• 42B37: "Harmonic analysis and PDEs" (MSC2020)

\end{abstract}

%%%%%%%%%%%%%%%%%%%%%%%%%%%%%%%%%%%%%%%%%%%%
%%%%%%%%%%%%%%%%%%%%%%%%%%%%%%%%%%%%%%%%%%%%
%%%%%%%%%%%%%%%%%%%%%%%%%%%%%%%%%%%%%%%%%%%%
%%%%%%%%%%%%%%%%%%%%%%%%%%%%%%%%%%%%%%%%%%%%
%%%%%%%%%%%%%%%%%%%%%%%%%%%%%%%%%%%%%%%%%%%%
%%%%%%%%%%%%%%%%%%%%%%%%%%%%%%%%%%%%%%%%%%%%
%%%%%%%%%%%%%%%%%%%%%%%%%%%%%%%%%%%%%%%%%%%%
%%%%%%%%%%%%%%%%%%%%%%%%%%%%%%%%%%%%%%%%%%%%
%%%%%%%%%%%%%%%%%%%%%%%%%%%%%%%%%%%%%%%%%%%%
%%%%%%%%%%%%%%%%%%%%%%%%%%%%%%%%%%%%%%%%%%%%

\section{Introduction}

\setcounter{theorem}{0}
\setcounter{equation}{0}

%%%%%%%%%%%%%%%%%%%%%%%%%%%%%%%%%%%%%%%%%%%%
%%%%%%%%%%%%%%% The Equation %%%%%%%%%%%%%%%
%%%%%%%%%%%%%%%%%%%%%%%%%%%%%%%%%%%%%%%%%%%%

\noindent
%Let $p\in(1,\infty]$.
We study solutions of the $p$-Laplace equation which yields
\begin{align}\label{eq:plapequation}
\Delta_{p} u :=\nabla \cdot ( |\nabla u |^{p - 2}\, \nabla  u ) = 0
\end{align}

\bigskip

\noindent
when $p \in(1, \infty)$.
If $p = \infty$, then the equation can be written
\begin{align}\label{eq:inflapequation}
\Delta_{\infty} u := \sum_{i,j = 1}^{n} \frac{\partial u}{\partial x_i} \frac{\partial u}{\partial x_j} \frac{\partial^2 u}{\partial x_i \partial x_j} = 0
\end{align}

\bigskip

\noindent
which is the so called $\infty$-Laplace equation.
%We refer the reader to Section \ref{sec:prel} for the definitions of weak solutions, viscosity solutions and $p$-harmonicity.
The $p$-Laplace equation arises in minimization problems,
nonlinear elasticity theory,
Hele-Shaw flows and
image processing, see e.g. %Lundstr\"om
\cite[Chapter 2]{avhandlingen} and the references therein for more on applications and motivations.

Let $\Omega\subset \mathbb{R}^n$ be a regular bounded domain and let $f$ be a real-valued continuous function defined on $\partial{\Omega}$.
It is well known that there exists a unique smooth function $u$, harmonic in $\Omega$,
such that $u = f$ continuously on $\partial \Omega$. The maximum principle and the Riesz representation theorem yield the following
representation formula for $u$,
\begin{equation*}
u(z)=\int\limits_{\partial\Omega}f(w) \,d{\omega}^{z}(w) \text{, \quad whenever\; $z\in \Omega$.}
\end{equation*}

\bigskip

\noindent
Here, $\omega^z(w) = \omega(dw, z, \Omega)$ is referred to as the harmonic measure at $z$ associated to the Laplace operator.
As the harmonic measure allows us to solve the Dirichlet problem, its properties are of fundamental
interest in classical potential theory.

In this paper we prove estimates in planar sectors of the following $p$-harmonic measure,  which is a generalization of harmonic measure, related to the $p$-Laplace equation.

%%%%%%%%%%%%%%%%%%%%%%%%%%%%%%%%%%%%%%%%%%%%
%%%%%%%%%%% p-harmonic measure %%%%%%%%%%%%%
%%%%%%%%%%%%%%%%%%%%%%%%%%%%%%%%%%%%%%%%%%%%

%
\begin{definition} \label{def:p-hmeas}
Let $G \subseteq \mathbb{R}^n$ be a domain, $E \subseteq \partial{G}$, $p \in (1,\infty)$
and $x \in G$.
The $p$-harmonic measure of $E$ at $x$ with respect to $G$, denoted by $\omega_p(x) = \omega_p(E,x,G)$,
is defined as $\inf_{u} u(x)$,
where the infimum is taken over all $p$-superharmonic functions $u \ge 0$ in $G$
such that $\liminf_{z \to y} u(z) \ge 1$, for all $y \in E$.
\end{definition}
\noindent
The $\infty$-harmonic measure is defined in a similar manner, but with $p$-superharmonicity replaced by absolutely minimizing
%, see Peres--Schramm--Sheffield--Wilson
\cite[pages 173--174]{inf-tow}.
It turns out that the $p$-harmonic measure in Definition \ref{def:p-hmeas} fails to be a measure but is a
$p$-harmonic function in $\Omega$, bounded below by $0$ and bounded above by $1$.
For these and other basic properties of $p$-harmonic measure we refer the reader to
%Heinonen--Kilpel\"ainen--Martio
\cite[Chapter 11]{HKM}.

Let $(r,\phi)$ be polar coordinates for $(x,y) \subset \mathbb{R}^2$
and consider the planar sector
\begin{align}\label{eq:def-S_v}
\mathcal{S}_\nu = \left\{ (x_1, x_2) \in  \mathbb{R}^2 \setminus \{(0, 0)\}; |\phi| <
\frac{\pi}{2\nu}\right\}, \quad \text{where} \quad \nu \geq \frac{1}{2},
\end{align}

\bigskip
\noindent
having aperture $\pi/\nu$ and apex at the origin.
%Put $\mathcal{S}_v(R) = B(0, R)\cap \mathcal{S}_v$ and
%$\Delta_v(R) = \partial B(0,R) \cap \mathcal{S}_v$.
Suppose that $p \in (1,\infty]$, $v \in [1/2,\infty)$ and
let $\omega_p(x) = \omega_p(\partial B(0,R) \cap \mathcal{S}_\nu, x,  B(0, R)\cap \mathcal{S}_\nu)$
%$\omega_p(\Delta_v(R), x,  \mathcal{S}_v(R))$
be the $p$-harmonic measure of $\partial B(0,R) \cap \mathcal{S}_\nu$ at $x$ with respect to $B(0, R)\cap \mathcal{S}_\nu$.
Using comparison arguments and basic boundary estimates together with certain explicit $p$-harmonic functions derived  in~\cite{A84}, \cite{A86}, \cite{stream} and \cite{persson-lic} we prove in Theorem \ref{th:p-harmonic-measure} that
\begin{align}\label{eq:estimate_1}
C^{-1} \left(\frac{|x|} {R}\right)^{k(\nu,p)} \,  \leq \omega_p(x)\, \leq C \left(\frac{|x|} {R}\right)^{k(\nu,p)},
\end{align}
whenever $x \in B(0,R)\cap\mathcal{S}_{2\nu}$ and where
$C$ depends only on $\nu$ and $p$.
As the $p$-Laplace equation is invariant under rotations, scaling and translations,
Theorem \ref{th:p-harmonic-measure} holds for any planar sector.
The exponent $k(\nu,p)$ is given by
\begin{equation}\label{eq:radialexponent}
k(\nu,p)=\frac{(\nu-1)\sqrt{(1-2\nu)(p-2)^2+\nu^2p^2}+(2-p)(1-2\nu)+\nu^2p}{2(p-1)(2\nu-1)},
\end{equation}

\bigskip
\noindent
interpreted as a limit when $\nu = 1/2$ and $p = \infty$ so that
%\komN{We should remark in the construction of explicit solutions that the case $v < 1$ is "new" (seems new but I am not sure)}

%
%$$
%k(v,\infty) = \frac{(v-1)|v-1| + v^2 + 2v -1}{2 (2v-1)}, \quad
%$$
%
\begin{equation}\label{eq:kinfhej}
k(1/2,p) = \frac{p-1}{p}  \quad \text{and} \quad
k(\nu,\infty) =\left\{
	\begin{array}{ll}
		1 \quad &\text{when} \quad \frac{1}{2} \leq \nu \leq 1,\\
		\frac{ \nu^2}{2\nu - 1} \quad  &\text{when} \quad 1 \leq \nu.
	\end{array}
\right.
\end{equation}
%
%$$
%k(v,\infty) = \frac{v^2}{2v - 1} \qquad \text{when $v \geq 1$ and 1 otherwise}
%$$

\begin{figure}[h]
\centering
\includegraphics[width=8cm,height=7cm]{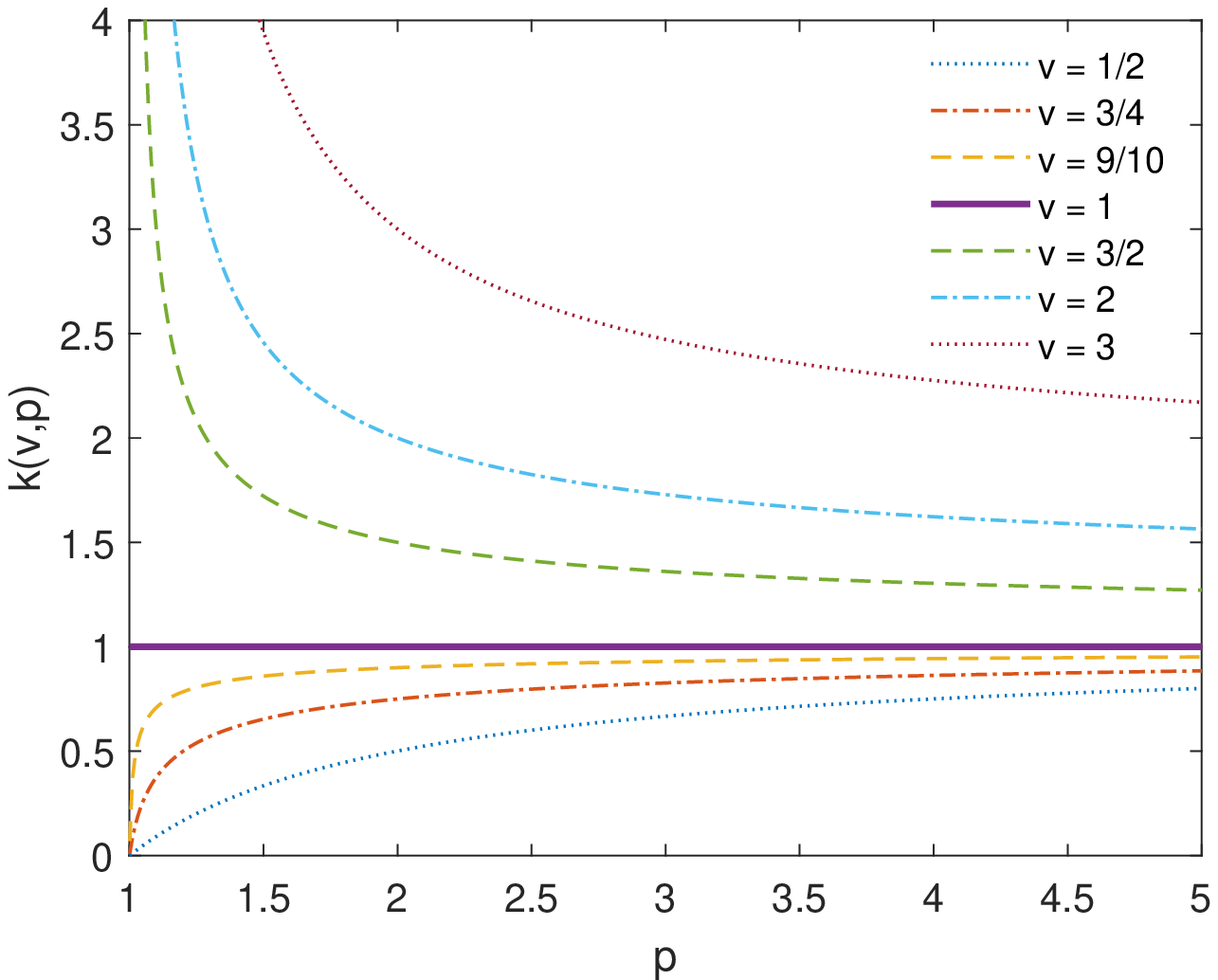}
\hspace{2mm}
\includegraphics[width=8cm,height=7cm]{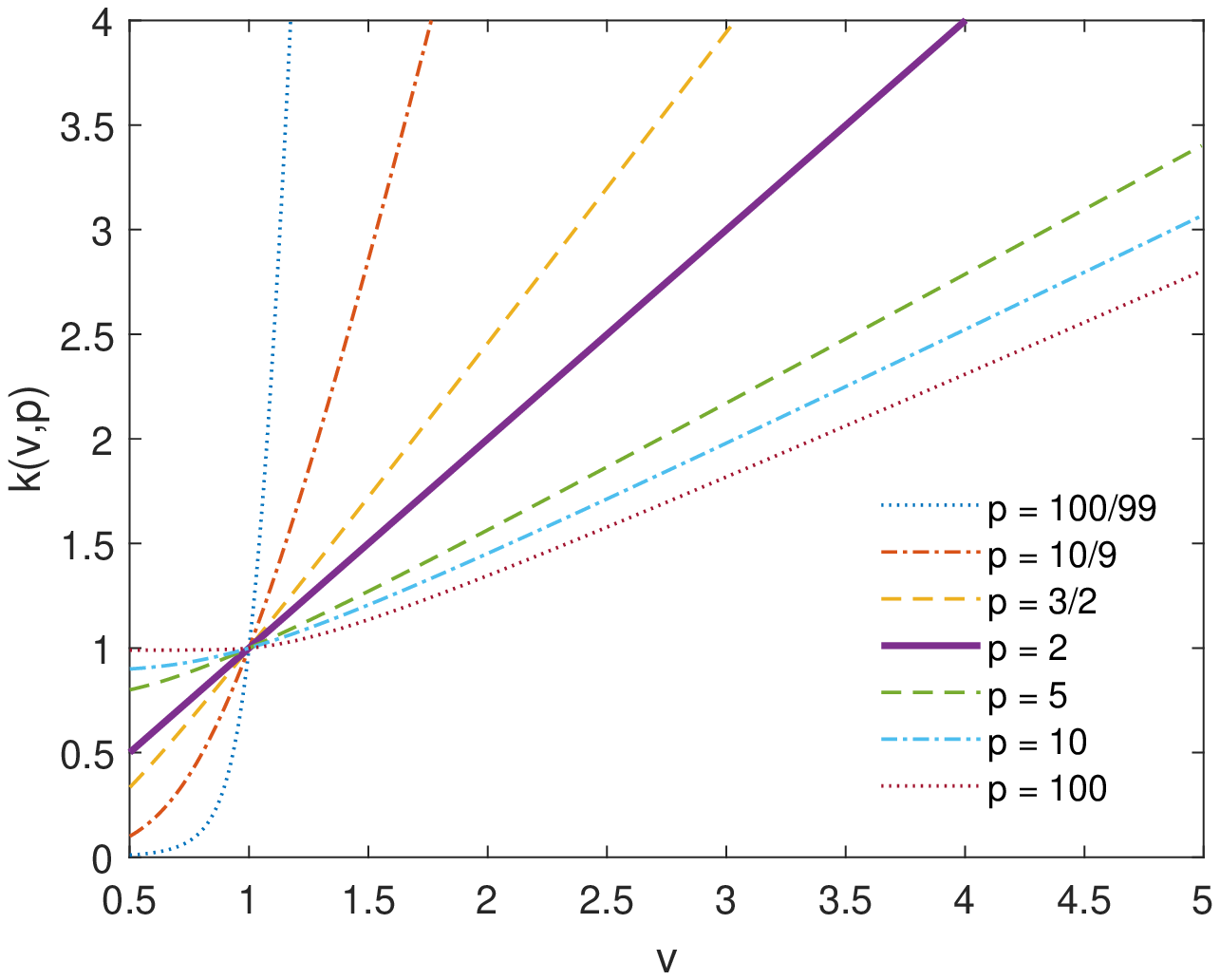}
\caption{The radial exponent $k(\nu,p)$ as function of $\nu$ and $p$.}
\label{fig:exponent_fig_1}
\end{figure}
Figure~\ref{fig:exponent_fig_1} shows the radial exponent $k(\nu,p)$ for some $\nu$ and $p$.
Curves for $\nu < 1$ approaches zero as $p\to 1$ and 1 as $p \to \infty$.
Curves for $\nu > 1$ approaches infinity as $p \to 1$ and $\nu^2/(2\nu - 1)$ as $p \to \infty$.
Moreover, $k(\nu,p) \to \infty$ as $\nu \to \infty$, reflecting the case when the sector $\mathcal{S}_\nu$ approaches a line.
The asymptotic behaviour in this case is $k(\nu,p) = \frac{p \nu}{2(p-1)} + O(1)$ which is in agreement with a related result in \cite{rysk}.
%
%We remark that the result $k(v,\infty) \equiv 1$ when $v \in [1/2,1]$ can also be proved using the comparison with cones principle for infinity harmonic functions, see Remark \ref{rem:infty}.
Further, the case $k(1/2,p)$ captures the rate at which the $p$-harmonic measure (or a positive $p$-harmonic function) vanish at a halfline because $\mathcal{S}_{1/2} = \mathbb{R}^2 \setminus\{(r,\phi) : \phi = \pi\}$.
Furthermore, we retrieve the known results in the classical cases
$
k(\nu,2) = \nu$ %\qquad\text{and}\qquad
and $k(1,p) = 1,
$
of which the first corresponds to the harmonic measure $(p = 2)$ and in the second  $\mathcal{S}_\nu$ is a halfplane. %(both cases should be well known).
%That $k(v,\infty) \equiv 1$ when $v \in [1/2,1]$ can be realized using the comparison with cones principle for infinity harmonic functions, see Remark \ref{rem_cones}.
%When $v=1$ then $\mathcal{S}_v$ is a halfplane and it should be well known that the exponent is $k(1,p) = 1$.
%

When $p = 2$ the $p$-harmonic measure coincides with the famous harmonic measure and our result, expressed in probabilistic terms, answers the question; what is the probability that a Brownian motion started at $x \in \mathcal{S}_\nu$ will first hit the part of the boundary consisting of the arc $\partial B(0,R)\cap \mathcal{S}_\nu$?
Our estimate in \eqref{eq:estimate_1} implies that the probability is comparable to $(|x|/R)^{\nu}$.

Our estimates for the $p$-harmonic measure imply local growth estimates %such as a boundary Harnack inequality
for $p$-sub- and $p$-superharmonic functions %(as well as for $p$-sub and $p$-superharmonic functions)
vanishing on a fraction of a domain contained in, or containing, a sector.
Indeed, we conclude that solutions must vanish at the same rate as
$|x|^{k(\nu,p)}$ as $x$ approaches the apex (Corollary \ref{cor:subsuper}).
Similar growth estimates %such as those in Corollary \ref{cor:subsuper}
where proved in the setting of $C^{1,1}$-domains in
%These estimates are similar to those proved for $C^{1,1}$-doamins in
%Aikawa--Kilpel\"ainen--Shanmugalingam--Zhong~
\cite{aikawa} and for wider classes of equations and other geometric settings in
%Lundstr\"om \cite{L11, L16}
%Lundstr\"om--Olofsson--Toivanen \cite{LOT20}.
\cite{LN07, LN08, LN10, L11, L16, AL16, LOT20}.
An immediate consequence of Corollary \ref{cor:subsuper} is the boundary Harnack inequality for $p$-harmonic functions in planar sectors, see \eqref{eq:second-cor}.
For $\nu \neq \frac{1}{2}$ such result is already well known by \cite{LN07,LN10} since then $\mathcal{S}_\nu$ is a Lipschitz domain.

Consider a domain $\Omega \subset \mathbb{R}^2$ having a sharp outwardly pointed cusp with apex $w$ and let $u$
be a $p$-subharmonic function taking nonpositive boundary values in a neighborhood of $w$.
Using Corollary \ref{cor:subsuper} we prove that then the rate of convergence to zero, as $x$ approaches the apex, is faster than any power of $|x-w|$,
i.e. for any $N > 0$ it holds that
%
%\begin{align*}
%\limsup_{\Omega \ni x\to w} \frac{u(x)}{|x-w|^{N}} \leq 0,
%\end{align*}
\begin{align*}
\underset{x\in\Omega}{\underset{x\to w}{\limsup}}\;\frac{u(x)}{|x-w|^N}\leq 0,
\end{align*}

\bigskip
\noindent
which is a result proved already in \cite[Theorem 3]{rysk}.
Consider now instead a domain $\Omega \subset \mathbb{R}^2$ having a sharp inwardly pointed cusp at $w$, and let $v$ be a $p$-superharmonic function taking nonnegative boundary values in a neighborhood of $w$.
In this case we prove that the rate of convergence to zero, as $x$ approaches the apex, is slower than $|x-w|^\frac{p-1}{p}$, i.e.,
%
%\begin{align*}
%\liminf_{\Lambda \ni x\to w} \frac{ v(x)}{|x-w|^\frac{p-1}{p}} > 0, \quad \text{where}\quad \Lambda = \{x \in \Omega : %d(x,\partial \Omega) \geq |x-w|\}.
%\end{align*}
%
\begin{align*}
\underset{x\in\Lambda}{\underset{x\to w}{\limsup}}\;\frac{u(x)}{|x-w|^N}>0, \quad \text{where}\quad \Lambda = \{x \in \Omega : d(x,\partial \Omega) \geq |x-w|\}.
\end{align*}

\bigskip

The $p$-harmonic measure has a probabilistic interpretation
in terms of the zero-sum two-player game \emph{tug-of-war},
see %Peres--Sheffield
\cite{noisy-tow} and \cite{inf-tow},
in which also estimates for $p$-harmonic measure are proved, e.g. for porous sets.
Further results in the literature include %Lundstr\"om--Vasilis
\cite{LV13} who proved estimates for $p$-harmonic measures in the plane,
which,
together with a result in %Hirata
\cite{H08},
yield properties of the $p$-Green function.
%
%DeBlassie--Smits
Estimates for the $p$-harmonic measure of a small spherical cap and of small axially symmetric sets are proved in \cite{DS16, DS18},
and in \cite{DS21} estimates for the $p$-harmonic measure is given of the part of the boundary of an infinite slab outside a cylinder.
In \cite{L16} estimates of $p$-harmonic measure, $n-m < p \leq \infty$, for sets in $\mathbb{R}^n$ which are close to an $m$-dimensional hyperplane, $0 \leq m \leq n-1$ are proved,
and in
%
%Llorente--Manfredi--Troy--Wu
\cite{LMTW19}
it is proved that
%,  for $1 < p < \infty$ and $N \geq 2$, by using ODE shooting techniques, that there is a constant $\alpha(p, N) > 0$ such that
the $p$-harmonic measure in  $\mathbb{R}_+^n$ of a ball of radius $0 < \delta \leq 1$ in $\mathbb{R}^{n-1}$ is bounded above and below by a constant times $\delta^{\alpha}$,
and explicit estimates for $\alpha$ are given.
For more on possible applications of $p$-harmonic measure, %has shown to be useful when estimating solutions to the $p$-Laplace equation,
see e.g. \cite[Chapter 11 and Chapter 14]{HKM}
%including the relation to Phragmen--Lindel\"ofs theorem which we will consider in this paper.
including Phragmen--Lindel\"ofs theorem and the study of quasiregular mappings.

%%%%%%%%%%%%%%%%%%%%%%%%%%%%%%%%%%%%%%%%%%%%
%%%%%%%% Phragmen-Lindelof Theorem %%%%%%%%%
%%%%%%%%%%%%%%%%%%%%%%%%%%%%%%%%%%%%%%%%%%%%

In Section \ref{sec:phragmen} we use the estimates in Theorem \ref{th:p-harmonic-measure} to prove Theorem \ref{th:phragmen}
which is an extended version of the
classical result of Phragm\'en–Lindel\"of \cite{PL08}.
In particular, suppose that $u$ is p-subharmonic in
an unbounded planar domain $\Omega$ contained in the sector $\mathcal{S}_\nu$
and suppose that $\limsup_{z\to\partial\Omega} u(z) \leq 0$.
Then either $u \leq 0$ in the whole of $\Omega$ or it holds that
\begin{align*}
\liminf_{R\to\infty} \left( \frac{1}{R^{k(\nu,p)}}  \sup_{\partial B(0,R) \cap \Omega}u\right) > 0,
\end{align*}

\bigskip

\noindent
where $k(\nu,p)$ is as in \eqref{eq:radialexponent}.
When $\Omega \equiv \mathcal{S}_\nu$,
the above growth rate is sharp.
We remark that when $\nu = 1$ the sector $\mathcal{S}_\nu$ is a halfplane and $k(1,p) = 1$; thus we retrieve the classical result that $p$-subharmonic functions must grow at least as fast as the distance to the boundary \cite{L85}.

In connection with the above Phragm\'en-Lindel\"of result we also prove,
for  $p \in (1,\infty)$, $\nu \in (1/2,\infty)$,
that positive $p$-harmonic functions in $\mathcal{S}_\nu$, vanishing on $\partial\mathcal{S}_\nu$, are unique (modulo normalization),
see Theorem \ref{th:unique}.
%This result is given in Theorem \ref{th:unique} and implies,
%that there exists a constant $C$ such that any such $p$-harmonic function
%can be written
%$
%u(x) = C r^k f_{v,p}(\phi)
%$.
The proof of this result uses %, besides the explicit $p$-harmonic functions \cite{A86},
scaling arguments and a boundary estimate from \cite{LN10}.

Being a generalization of maximum principles to unbounded domains
the Phragm\'en-Lindel\"of principle \cite{PL08} is undoubtedly an important result
%have been frequently studied during the last century.
%One reason may be its
with applications in e.g. elasticity theory \cite{Horgan, Q93,Qnew}.
To summarize some literature (without giving a complete list)
we mention that results of \cite{PL08} was extended to halfspaces of
$\mathbb{R}^n$ in \cite{A37} and to
 general elliptic equations of second order in \cite{G52, S54,H64}.
Uniformly elliptic equations in nondivergence form in cones were considered in \cite{M71},
growth estimates of bounded solutions of quasilinear equations in \cite{K93,JL03} and for elliptic equations in sectors, see \cite{V04}.
%Capuzzo-Dolcetta--Vitolo \cite{CDV07} and %Armstrong--Sirakov--Smart \cite{ASS12} considered
Fully nonlinear equations were considered in
\cite{CDV07,ASS12}, the later in certain Lipschitz domains, and
 \cite{KN09} considered fully nonlinear elliptic PDEs with unbounded coefficients and nonhomogeneous terms.
%Adamowicz \cite{A14} studied subsolutions of
Results for variable exponent $p$-Laplace-type equations can be found in \cite{A14},
while %Bhattacharya \cite{Bhatt05} and Granlund--Marola \cite{GM14} considered
infinity-harmonic functions are considered in \cite{Bhatt05, GM14}.
In \cite{L85}, Phragm\'en-Lindel\"of's theorem for $n$-subharmonic functions, when the boundary is an $m$-dimensional hyperplane in $\mathbb{R}^n$,
$0 \leq m \leq n-1$ is proved.
This was extended to $p$-subharmonic functions,
$n-m < p \leq \infty$, in \cite{L16}.
In \cite{BM20} it is showed that solutions of a generalized $p$-Laplace equation in the upper halfplane, vanishing on $\{x_n = 0\}$, equals $u(x)= x_n$ (modulo normalization), while
 the growth of solutions of the minimal surface equation over domains containing a halfplane was considered in \cite{LW15}.
A Phragm\'en-Lindel\"of theorem for a mixed boundary value problem for  quasilinear elliptic equations of $p$-Laplace type in an open infinite circular half-cylinder was proved in \cite{BM21}.
The spatial behavior of solutions of
the Laplace equation on a semi-infinite cylinder with dynamical nonlinear
boundary conditions was investigated in \cite{Qnew}.
In halfspaces of $\mathbb{R}^n$, growth estimates for subsolutions of fully nonlinear nonhomogeneous PDEs was characterized in terms of solutions to certain ordinary differential equations in \cite{L21}.
Using this characterization, several
Phragmen-Lindel\"of-type results were derived, e.g. a sharp theorem for the variable exponent $p$-Laplace equation, and also sharp results for equations allowing for sublinear growth in the gradient.
Phragm\'en-Lindel\"of theorems for plurisubharmonic functions on cones were proved in \cite{M92} while $k$-Hessian equations with lower order terms were considered in \cite{BhM21}.
The present paper complements the above by giving the sharp exponent $k(\nu,p)$ explicitly in case of positive $p$-harmonic functions in planar sectors.

In Section \ref{sec:prel} we summarize some well known basic definitions and properties of solutions to the $p$-Laplace equation and  in Sections \ref{sec:explicit} and \ref{sec:appendix} we summarize and prove properties on explicit $p$-harmonic functions in planar sectors.
Using these results we state and prove our estimates of $p$-harmonic measure in Section \ref{sec:est-p-meas}, while in Section \ref{sec:est-p-harmonic} we give Corollaries for $p$-sub- and $p$-superharmonic functions.
Sections \ref{sec:phragmen} and \ref{sec:unique} is devoted to Phragmen-Lindel\"ofs theorem and uniqueness of $p$-harmonic functions in sectors, respectively.
We end the paper by showing in Section \ref{sec:inf} that in the case of infinity-harmonic measure and infinity-harmonic functions %(i.e. $p = \infty$)
then most of our results extend to symmetric $n$-dimensional domains, $n > 2$.

%%%%%%%%%%%%%%%%%%%%%%%%%%%%%%%%%%%%%%%%%%%%%%%%%%%%%%%%%%%%%%%%%%%%%%%%%%%%%%%%%%%%
%%%%%%%%%%%%%%%%%%%%%%%%%%%%%%%%%%%%%%%%%%%%%%%%%%%%%%%%%%%%%%%%%%%%%%%%%%%%%%%%%%%%
%%%%%%%%%%%%%%%%%%%%%%%%%%%%%%%%%%%%%%%%%%%%%%%%%%%%%%%%%%%%%%%%%%%%%%%%%%%%%%%%%%%%
%%%%%%%%%%%%%%%%%%%%%%%%%%%%%%%%%%%%%%%%%%%%%%%%%%%%%%%%%%%%%%%%%%%%%%%%%%%%%%%%%%%%
%%%%%%%%%%%%%%%%%%%%%%%%%%%%%%%%%%%%%%%%%%%%%%%%%%%%%%%%%%%%%%%%%%%%%%%%%%%%%%%%%%%%
%%%%%%%%%%%%%%%%%%%%%%%%%%%%%%%%%%%%%%%%%%%%%%%%%%%%%%%%%%%%%%%%%%%%%%%%%%%%%%%%%%%%
%%%%%%%%%%%%%%%%%%%%%%%%%%%%%%%%%%%%%%%%%%%%%%%%%%%%%%%%%%%%%%%%%%%%%%%%%%%%%%%%%%%%
%%%%%%%%%%%%%%%%%%%%%%%%%%%%%%%%%%%%%%%%%%%%%%%%%%%%%%%%%%%%%%%%%%%%%%%%%%%%%%%%%%%%
%%%%%%%%%%%%%%%%%%%%%%%%%%%%%%%%%%%%%%%%%%%%%%%%%%%%%%%%%%%%%%%%%%%%%%%%%%%%%%%%%%%%
%%%%%%%%%%%%%%%%%%%%%%%%%%%%%%%%%%%%%%%%%%%%%%%%%%%%%%%%%%%%%%%%%%%%%%%%%%%%%%%%%%%%

\section{Preliminaries}
\label{sec:prel}

\setcounter{theorem}{0}
\setcounter{equation}{0}

\noindent
In this section we state some basic definitions and results for $p$-harmonic measure and $p$-harmonic functions needed later.
By $\Omega$ we denote a domain, that is, an open connected set.
For a set $E \subset \mathbb{R}^n$
we let $\overline{E}$ denote the closure and
$\partial E$ the boundary. % and $\complement E$ the complement of $E$. %and we put $E^o = E\setminus \partial E$.
Further, $d(x,E)$ denotes the Euclidean distance from $x \in \mathbb{R}^n$ to $E$, and
$B(z,r) = \{ x\in\mathbb{R}^2 : d(x,z) < r \}$
denotes the open ball with radius $r$ and center $z$.
By $c$ we denote a constant $\geq 1$, not necessarily the same at each occurrence,
depending only on $\nu$ and $p$ if nothing else is mentioned.
Moreover, $c(a_1, a_2, \dots, a_k)$ denotes a constant $\geq 1$, not necessarily the same at each occurrence,
depending only on $a_1, a_2, \dots, a_k$.
%and we write $A \approx B$ if there exists a constant $c$ such that $c^{-1}A \leq B\leq cA$.
%We denote points in Euclidean $n$-space $\mathbf R^{n}$ by
%$x = (x_1,x_2,\dots,x_n) = (x',x'')$, where
%%
%\begin{align}\label{eq:notation-x'-and-x''}
%x' = (x_1, x_2, \dots, x_{n-m}) \quad \textrm{and} \quad x'' = (x_{n-m+1}, x_{n-m+2},\dots, x_n).
%\end{align}
%%
%Finally, we write $N = \{1,2,3,\dots\}$ for the set of natural numbers.

We proceed with defining weak and viscosity  solutions and $p$-harmonicity.
If $p \in (1, \infty)$, we say that $u$ is a \textit{(weak) subsolution (supersolution)} to the $p$-Laplacian in
a domain $\Omega$ provided $u \in W_{loc}^{1,p}(\Omega)$ and
\begin{equation*}
\int\limits_{\Omega} | \nabla  u |^{p - 2}  \, \langle   \nabla  u , \nabla \theta  \rangle \, dx \leq (\geq) \,0,
\end{equation*}

\bigskip

\noindent
whenever $\theta \in C^{\infty}_0(\Omega)$ is non-negative.
A function $u$ is a \textit{(weak) solution} of the $p$-Laplacian if it is both a subsolution and a supersolution.
Here, and in the sequel,
$W^{1,p}(\Omega)$ is the Sobolev space of those $p$-integrable functions whose first
distributional derivatives are also $p$-integrable,
%$W^{1,q}_0(\Omega)$ denotes the set of functions in $W^{1,p}(\Omega)$ with compact support in $\Omega$.
and $C^\infty_0(\Omega)$ is the set of infinitely differentiable functions with compact support in $\Omega$.
If $p = \infty$, the equation is no longer of divergence form and therefore the above definition needs to be replaced.
We use instead the notion of viscosity solutions.
Here, and in the sequel, $\Delta_{\infty}$ is the $\infty$-Laplace operator defined in \eqref{eq:inflapequation}.

An upper semicontinuous function $u : \Omega \rightarrow \mathbb{R}$ is a \textit{(viscosity) subsolution} of the
$\infty$-Laplacian in $\Omega$ provided that for each function $\psi \in C^2(\Omega)$ such that $u-\psi$ has a local maximum at a point $x_0 \in \Omega$, we have $\Delta_{\infty}\psi(x_0) \geq 0$.
A lower semicontinuous function $u : \Omega \rightarrow \mathbb{R}$ is a \textit{(viscosity) supersolution} of the
$\infty$-Laplacian in $\Omega$ provided that for each function $\psi \in C^2(\Omega)$ such that $u-\psi$ has a local minimum at a point $x_0 \in \Omega$, we have $\Delta_{\infty}\psi(x_0) \leq 0$.
A function $u : \Omega \rightarrow \mathbb{R}$ is a \textit{(viscosity) solution} of the
$\infty$-Laplacian if it is both a subsolution and a supersolution.

If $u$ is an upper semicontinuous subsolution to the $p$-Laplacian in $\Omega$, $p \in (1, \infty]$,
then we say that $u$ is \textit{$p$-subharmonic} in $\Omega$.
If $u$ is a lower semicontinuous supersolution to the $p$-Laplacian in $\Omega$, $p \in (1, \infty]$,
then we say that $u$ is \textit{$p$-superharmonic} in $\Omega$.
If $u$ is a continuous solution to the $p$-Laplacian in $\Omega$, $p \in (1, \infty]$, then $u$ is
\textit{$p$-harmonic} in $\Omega$.

We note that for the $p$-Laplacian, $1 < p < \infty$, weak solutions are also viscosity solutions
(defined as above but with $\Delta_{\infty}$ replaced by $\Delta_p$); see \cite[Theorem 1.29]{Ju}.
Moreover, under suitable assumptions, an $\infty$-harmonic function is the uniform limit of a sequence of $p$-harmonic functions as $p \to \infty$; see \cite{J}.
For more on weak solutions, viscosity solutions, $p$-harmonicity and $p$-superharmonicity, see for instance \cite{HKM} and \cite{guide}.

We will make use of the nowadays well known basic properties of $p$-harmonic functions:

\begin{lemma} \label{jamforelseprin}
Let $p \in (1,\infty]$ and suppose that $u$ is $p$-superharmonic and that $v$ is $p$-subharmonic in a bounded domain  $\Omega \subset \mathbb{R}^n$. If
\begin{equation*}
\limsup_{z\to w} v(z) \leq \liminf_{z\to w} u(z)
\end{equation*}
for all $w \in \partial\Omega$,
and if both sides of the above inequality are not simultaneously $\infty$ or $-\infty$,
then $v \leq u$ in $\Omega$.
\end{lemma}
\noindent
\begin{proof}
%{\bf Proof.}
If $p \in (1, \infty)$ then this follows from \cite[Theorem~7.6]{HKM}.
For the case $p = \infty$ this was first proved in \cite[Theorem 3.11]{J}.
A shorter proof was later given in \cite{armstrongsmart}.\end{proof}

\begin{lemma}\label{le:scaling}
Let $\Omega \subset \mathbb{R}^n$ be a domain, $p \in (1,\infty]$ and assume
that $u$ is a $p$-harmonic function in $\Omega$.
Let $k \in \mathbb{R}$ and $z \in \mathbb{R}^n$.
Then $\hat u$ is $p$-harmonic in some $\widetilde \Omega \in \mathbb{R}^n$,
in either of the following cases:
$$
(i)\; \hat u = ku, \quad (ii)\; \hat u = u(x + z),\quad (iii)\; \hat u = u(kx).
$$
\end{lemma}

\noindent
\begin{proof}
Follows by standard calculations.\end{proof}

\begin{lemma}\label{harnack}
 Let $p \in (1, \infty]$, $w \in \mathbb{R}^n$, $r\in (0,\infty)$ and suppose that $u$ is a positive $p$-harmonic function in $B(w,2r)$. Then there exists a constant $c \in (1,\infty]$, depending only on $p$ and $n$, such that
\begin{equation*}
\sup_{B(w, r)} u \leq c \inf_{B(w, r)} u \text{.}
\end{equation*}
\end{lemma}

\noindent
\begin{proof}
%{\bf Proof.}
For the case $p \in (1,\infty)$, see e.g. \cite{KMV}. For the case $p = \infty$ the result follows by taking the limit $p\to \infty$ in the former case, see \cite{lin-man}.%, or \cite[Lemma~2.2]{LN}.
\end{proof}

\bigskip

The following well known estimate tells that $p$-harmonic functions are H\"older continuous up to the boundary in the setting of the rather general class of non-tangentially accessible (NTA) domains.
We will only apply the result in smooth planar domains,
and refer the interested reader to e.g.~\cite[Chapter 1.6]{avhandlingen} for a definition of NTA-domains.

\begin{lemma}\label{le:holder}
Assume that $\Omega \subset \mathbb{R}^n$ is an NTA-domain with constant $M$,
let $w\in \partial \Omega$,
$0 <r< r_0$
and suppose that $u$ is a positive $p$-harmonic function
in $\Omega\cap B(w,2r)$,
continuous on $\overline{\Omega}\cap B(w,2r)$ with $u = 0$ on $\partial \Omega \cap B(w,2r)$.
Then there exist $c$ and $\alpha\in (0, 1]$,
depending only on $M$, $n$ and $p$, such that
\begin{align*}
|u(x) - u(y)| \leq c \left( \frac{|x - y|}{r} \right)^{\alpha} \sup_{B(w, 2r)\cap\Omega} u,
\end{align*}

\bigskip
\noindent
whenever $x,y \in B(w,r) \cap \Omega$.
\end{lemma}

\noindent
\begin{proof}
By observing that $\Omega$ is $p$-regular
by the NTA-assumption, the lemma follows by the same arguments as
in \cite[Theorem 6.44, Lemma 6.47]{HKM}.\end{proof}

\bigskip

The following Lemma states that any positive $p$-harmonic function, vanishing on a portion of the boundary of a $C^{1,1}$-domain,
must vanish at the same rate as the distance to the boundary.
The right inequality in Lemma \ref{le:bhi-C11} -- which is an immediate consequence of the left inequality -- is usually referred to as a boundary Harnack inequality and states that any two $p$-harmonic functions, vanishing on the boundary, must vanish at the same rate.

\begin{lemma}\label{le:bhi-C11}
Let $\Omega \in \mathbb{R}^n$ be a $C^{1,1}$-domain, or equivalently a domain satisfying the ball condition with radius $r_b$, $p \in (1,\infty]$, $n\geq 2$, $w \in \partial \Omega$ and $0 < r < r_b$.
Suppose that $u$ and $v$ are positive $p$-harmonic functions in $\Omega \cap B(w, r)$,
satisfying $u = 0 = v$ on $\partial \Omega \cap B(w, r)$.
Then there exists $c = c(n,p)$ such that
\begin{align*}
 c^{-1} \frac{d(x,\partial \Omega)}{r} \leq \frac{u(x)}{u(a_r(w))} \leq c  \frac{d(x,\partial \Omega)}{r}  \quad \text{and} \quad c^{-1} \frac{u(a_r(w))}{v(a_r(w))} \leq \frac{u(x)}{v(x)} \leq c  \frac{u(a_r(w))}{v(a_r(w))},
\end{align*}

\bigskip
\noindent
whenever $x \in \Omega \cap B(w, r/c)$.
Here, $a_r(w)$ is a point in $\Omega$
%guaranteed to exist in $C^{1,1}$-domains,
satisfying $d(a_r(w),w) = r/c$ and $d(a_r(w),\partial \Omega) = r/c$.
\end{lemma}

\noindent
\begin{proof}
For $p \in (1,\infty)$ we refer to \cite[Lemma 3.1 and Lemma 3.3]{aikawa}. See also \cite{LN07,LN08, LN10} for similar as well as stronger estimates in more general geometries.
When $p = \infty$ the proof is similar, but then the comparison argument uses cones (which are $\infty$-harmonic) in place of the function $\phi_p(x)$ defined on \cite[page 286]{aikawa}, and therefore the exterior ball condition is not needed. See also %Bhattacharya
\cite{Bhatt07} for the case $p=\infty$.
\end{proof}

\section{Explicit $p$-harmonic functions in sectors}
\label{sec:explicit}

\setcounter{theorem}{0}
\setcounter{equation}{0}

In this section we are gonna prove the following Lemma,
which is similar to \cite[Lemma 4.1]{LV13}, using some explicit $p$-harmonic functions derived in~\cite{A84}, \cite{A86}, \cite{stream} and \cite{persson-lic}. The result will be of crucial importance when we prove our main results in Sections~\ref{sec:est-p-meas}-\ref{sec:inf}.

\begin{lemma}\label{derivatalemma}
Let $p\in(1,\infty]$ and $\nu\in[\frac{1}{2},\infty)$.
Then there exists a positive $p$-harmonic function $u_{\nu,p}:\mathcal{S}_\nu\rightarrow \mathbb{R}$
of the form $u_{\nu,p}(x)=r^kf_{\nu,p}(\phi)$, where the exponent $k = k(\nu,p)$ is given by \eqref{eq:radialexponent}.
The exponent $k$ is non-decreasing in $\nu$,
increasing in $p$ for $\nu\in [1/2,1)$, and decreasing in $p$ for $\nu > 1$.
Moreover, the function $f_{\nu,p}(\phi)$ is continuous, differentiable and satisfies
\begin{itemize}\itemsep2mm
\item[$i)$] $f_{\nu,p}(\pm\frac{\pi}{2\nu})=f'_{\nu,p}(0) = 0$,
\item[$ii)$] $0 \leq f_{\nu,p}(\phi) \leq 1,\; |f'_{\nu,p}(\phi)|\leq c$ when $\phi\in[-\frac{\pi}{2\nu},\frac{\pi}{2\nu}]$,
\item[$iii)$] $|f_{\nu,p}(\phi)| \geq c^{-1}$ when  $\phi\in[-\frac{\pi}{4\nu},\frac{\pi}{4\nu}]$ and $|f'_{\nu,p}(\phi)| \geq c^{-1}$ when  $\phi\in[-\frac{\pi}{2\nu},\frac{\pi}{2\nu}] \setminus[-\frac{\pi}{4\nu},\frac{\pi}{4\nu}]$.
%$\exists\, m\in\mathbb{R}: m\leq f'(\theta)$, and %$f'_{\nu,p}(\pm\frac{\pi}{2v})\neq 0$.
%\item[$iii)$] \komN{State that derivative is bounded below near the boundary, as we need in proof of Theorem \ref{th:p-harmonic-measure}.}
\end{itemize}
\end{lemma}
%Harnacks
%\begin{equation}
%2^{nN}u(y)\geq2^{nN}\inf_{V}u(y)\geq %u(x)\geq\frac{1}{2^{nN}}\sup_{V}u(y)\geq\frac{1}{2^{nN}}u(y)
%\end{equation}
%
%We will prove Lemma~\ref{derivatalemma}, which is similar to \cite[Lemma 4.1]{LV13}, using some explicit $p$-harmonic functions derived in \cite{A86} and \cite{Pe}.
%The proof is split into four different cases; $p = 2$, $2 < p$, $p = \infty$ and $1 < p < 2$.
%
%We only present a short proof here,
%additional details and derivations are given in the Appendix \ref{sec:appendix}.  \\

\bigskip

\noindent
\begin{proof} % of Lemma~\ref{derivatalemma}.}
We begin by noting that standard calculations yield (see Appendix~\ref{find}) $dk/d\nu \geq 0$ whenever $\nu \in [\frac{1}{2},\infty)$ and $p\in(1,\infty]$.
Moreover, $dk/dp > 0$ for $\nu\in [1/2,1)$, and $dk/dp < 0$ for $\nu > 1$, whenever $p\in(1,\infty)$.
The rest of the proof is split into four different cases; $p = 2$, $2 < p$, $p = \infty$ and $1 < p < 2$.

\subsubsection*{Case 1: $p=2$.}

In this case \eqref{eq:plapequation}
reduces to the Laplace equation which yields, in polar coordinates, %$x(r,\phi)=r \cos \phi$, $y(r,\phi)=r\sin\phi$,
\begin{align*}
   % \frac{\partial^2 u}{\partial r^2} + \frac{1}{r} \frac{\partial u}{\partial r} + \frac{1}{r^2}\frac{\partial^2 u}{\partial \phi^2} = 0.
    u_{rr} + \frac{1}{r} u_r + \frac{1}{r^2} u_{\phi \phi} = 0.
\end{align*}

\bigskip
\noindent
Thus, the solution $u(r,\phi)=r^\nu\cos(\nu\phi)$ has the desired properties since $k(\nu,p) = \nu$.
%Below we will also see that this solution coincides with the more general general solutions in the cases $1<p<2$ and $2<p<\infty$ when we approaches $p=2$ in limiting processes.

%%%%%%%%%%%%%%%%%%%%%%%%%%%%%%%%%%%%%%%%%%%%%%%%%%%%%%%%%%%%%%%%%%%%%%%%%%%
%%%%%%%%%%%%%%%%%%%%%%%%%%%%%%%%%%%%%%%%%%%%%%%%%%%%%%%%%%%%%%%%%%%%%%%%%%%
%%%%%%%%%%%%%%%%%%%%%%%%%%%%%%%%%%%%%%%%%%%%%%%%%%%%%%%%%%%%%%%%%%%%%%%%%%%

\subsubsection*{Case 2: $2 < p < \infty$.}
In polar coordinates the $p$-Laplace equation, with $b = 1/(p-2)$, $p \neq 2$, yields (see Appendix~\ref{transform} for a derivation),
\begin{multline}\label{eq:polar-coords}
\left(b+1\right)u_r^2u_{rr}+\frac{b}{r^2}\left(u_{rr}u_{\phi}^2+u_r^2u_{\phi\phi}\right)+ \frac{(b+1)}{r^4}u_{\phi}^2u_{\phi\phi}+\frac{b}{r}u_r^3+\frac{(b-1)}{r^3}u_ru_{\phi}^2+\frac{2}{r^2}u_ru_\phi u_{r\phi}=0.
\end{multline}
We are searching for solutions of the form $u(r,\phi)=r^kf(\phi)$, where  $f(\phi)\in \mathcal{C}^2$ and $k$ are to be determined.
Inserted in \ref{eq:polar-coords} we obtain
\begin{equation}\label{eq:sepeqinsec3}
\left[\left(b+1\right)(f')^2+bk^2f^2\right]f''+\left(2k+bk-1\right)kf(f')^2+\left(bk+k-1\right)k^3f^3=0.
\end{equation}

\bigskip
\noindent
Equation~(\ref{eq:sepeqinsec3}) can be solved (for details check~\cite[Lemma 2]{A86}) to yield
\noindent
\begin{equation}\label{eq:losningv}
f_{\nu,p}(\phi) = c\left(1-\frac{\cos^2\theta_{\nu,p}(\phi)}{ak}\right)^{\frac{k-1}{2}}\cos\theta_{\nu,p}(\phi),
\end{equation}
where $a = (p-1)/(p-2)$ and $\theta_{\nu,p}$ is a certain continuous, strictly increasing function of $\phi$ and $c = c(\nu,p)$ is chosen so that $0 \leq f_{\nu,p} \leq 1$.
When $|\phi|<\frac{\pi}{2\nu}$, we have
(See Appendix~\ref{find}),
\begin{equation}\label{eq:implicitangel}
  \phi=\theta_{\nu,p}(\phi)-\left(1-\frac{1}{k}\right)\frac{\sqrt{ak}}{\sqrt{ak-1}}\left[\arctan\left(\lambda_{\nu,p}\tan\frac{\theta_{\nu,p}(\phi)}{2}\right)+\arctan\left(\frac{1}{\lambda_{\nu,p}}\tan\frac{\theta_{\nu,p}(\phi)}{2}\right)\right],
\end{equation}

\bigskip
\noindent
where $\lambda_{\nu,p} = \frac{\sqrt{ak - 1}}{\sqrt{ak} + 1}$.
The function $\theta_{\nu,p}$ is chosen so that it maps the interval $[\frac{-\pi}{2\nu}, \frac{\pi}{2\nu}]$ to $[-\frac{\pi}{2},\frac{\pi}{2}]$, and this condition determines the radial exponent $k$.
More precisely, the condition that determines $k$ is given by %\cite[X]{A86}, and reads
\begin{equation*}
  \frac{\pi}{\nu}=\phi\left(\frac{\pi}{2}\right)-\phi\left(-\frac{\pi}{2}\right)=\pi\left(1-\left(1-\frac{1}{k}\right)\frac{\sqrt{ak}}{\sqrt{ak-1}}\right).
\end{equation*}
%
%\bigskip
%
Solving for $k = k(\nu,p)$ we obtain the exponent given by \eqref{eq:radialexponent} in the introduction.
Moreover, $f_{\nu,p}\left(\pm\frac{\pi}{2\nu}\right) = 0$ and $ak > 1$ implying $f_{\nu,p}(\phi) > 0$ for $|\phi| < \frac{\pi}{2\nu}$.

\bigskip

We also need to estimate the derivative of $f_{\nu,p}$.
Differentiation of $f_{\nu,p}$ in \eqref{eq:losningv} and simplifying (see~\cite{A86} page 143, and/or Appendix~\ref{compute}), yield

\begin{align}\label{eq:bounded derivative}
    f_{\nu,p}'(\phi)=-k\, c \left(1-\frac{\cos^2\theta_{\nu,p}(\phi)}{ak}\right)^{\frac{k-1}{2}}\sin\theta_{\nu,p}(\phi).
\end{align}

\bigskip
\noindent
Since $ak > 1$ we have that $|f_{\nu,p}'(\phi)|\leq c(\nu,p)$ for all $\phi\in[-\frac{\pi}{2\nu},\frac{\pi}{2\nu}]$. Also $f_{\nu,p}'(\phi)=0$ will only occur when $\theta=n\pi$, $n\in \mathbb{Z}$, corresponding to $\phi=\frac{n\pi}{\nu}$.
Hence the only place where $f_{\nu,p}'$ is zero in $\mathcal S_\nu$ is in the radial direction along $\phi=0$.
It follows that we can conclude, by continuity of $f'_{\nu,p}$ which holds since $ak >1$, the existence of a constant $c = c(\nu,p)$ such that $|f'_{\nu,p}(\phi)| \geq c^{-1}$ whenever  $\phi\in[-\frac{\pi}{2\nu},\frac{\pi}{2\nu}] \setminus[-\frac{\pi}{4\nu},\frac{\pi}{4\nu}]$.
It also follows that $f_{\nu,p}(\phi) \geq c^{-1}$ whenever $\phi \in [-\frac{\pi}{4\nu},\frac{\pi}{4\nu}]$.
The proof when $2 < p<\infty$ is complete.

\subsubsection*{Case 3: $p=\infty$.}

Letting $a=1$ when $p=\infty$, the function from Case 2 is immediately extended to the case when $p=\infty$.
Indeed, in this case the separation equation~(\ref{eq:sepeqinsec3})  boils down to
\begin{equation}\label{eq:sep-inf}
(f')^2f''+\left(2k-1\right) kf(f')^2+\left(k-1\right)k^3f^3=0,
\end{equation}
%
%and the parametric representation can be written
which has solution
\begin{equation*}
%\left\{
%	\begin{array}{ll}
%		\phi=\theta^{\ast}+\displaystyle \int_{\theta^{\ast}}^\theta\frac{\sin^2\theta'}{k-\cos^2\theta'}\,d\theta'  \\[24pt]
		f_{\nu,\infty}(\phi)=c\left(1-\frac{\cos^2\theta(\phi)}{k}\right)^\frac{k-1}{2}\cos\theta(\phi)
%	\end{array}.
%\right.
\end{equation*}
with radial exponent
\begin{equation*}
k(\nu,\infty) =\left\{
	\begin{array}{ll}
		1 \quad &\text{when} \quad \frac{1}{2} \leq \nu \leq 1,\\
		\frac{ \nu^2}{2\nu - 1} \quad  &\text{when} \quad 1 \leq \nu.
	\end{array}
\right.
\end{equation*}

\bigskip

This case is studied in detail in~\cite{A84} and
$f_{\nu,\infty}$ is infinitely differentiable on $[-\frac{\pi}{2\nu},0) \cup(0,\frac{\pi}{2\nu}]$ and differentiable at $\phi = 0$.
%$f\in C^1[0,\bar{\phi }}]\,\cap\, C^{\infty}(0,\bar{\phi}})$, \komN{It is smooth except at points where $f' = 0$}
%where $\bar{\phi}}=\pi\left(1-\sqrt{1-\frac{1}{k}\right)$, since $a=1$.
It follows that $u_{\nu,\infty} = r^k f_{\nu,\infty}$ satisfies the required conditions, except that it is not immediate that the function is $\infty$-harmonic in the viscosity sense in the entire sector.
This is because at points where $f_{\nu,\infty}'(\phi) = 0$ we have $f''_{\nu,\infty}(\phi) = -\infty$, which can be seen from \eqref{eq:sep-inf}, and $f_{\nu,\infty}(\phi)$ is thus not $C^2$ there.
For $\nu = 1$, it was shown in \cite[appendix~I]{bhatta04} that $r^{-1/3}f_{1,\infty}$
is indeed $\infty$-harmonic in the viscosity sense. However, the exact same
proof works for all $\nu \ge \frac{1}{2}$; or more generally, for $\infty$-harmonic functions
with polar representation $r^{k} f(\phi) \ge 0$, with $k \cdot (1-k) \leq 0$
and where $f \in C^{1}$ is $C^2$ for $\phi \ne 0$, has a local maximum at $\phi = 0$ and
satisfies $\lim_{\phi \to 0} f''(\phi) = -\infty$.
To complete the proof of the lemma when $p = \infty$ we observe that statements $i)-iii)$ follow in a similar way as in the case $2 < p < \infty$.

\subsubsection*{Case 4: $1<p<2$.}

We are going via a \emph{stream function} technique to handle this situation.
Indeed, we will use the stream function for the case $2<p<\infty$ to find our desired solution for $1<p<2$.
We begin by the following lemma from~\cite{stream} for which we present a proof in Appendix~\ref{taste}.
\begin{lemma}\label{le:konjugatlemma}
 Let $u(r,\phi)=r^kf(\phi)$ be $p$-harmonic in sector $\mathcal{S}_\nu$, $k>0$, and $2<p<\infty$. Then there exists a $q$-harmonic stream function $v(r,\phi)=r^{\lambda}g(\phi)$ in  $\mathcal{S}_\nu$, where $\lambda=(p-1)(k-1)+1$, $q=p/(p-1)$, and
\begin{equation*}%\label{eq:streamphi}
     g(\phi)=-\dfrac{1}{\lambda}f'(\phi)\left(k^2f(\phi)^2+f'(\phi)^2\right)^{\frac{p-2}{2}}.
\end{equation*}

\bigskip
\noindent
The function $g(\phi)$ is periodic whenever $f(\phi)$ is.
\end{lemma}
To apply Lemma \ref{le:konjugatlemma} it is convenient to view the $p$-values at hand as conjugate to those in Case 2, i.e. when $2 < p < \infty$.
Hence, fix $q \in (1, 2)$ and denote the conjugate index by $p = q / (q-1) \in (2,\infty)$.
Substituting $p=q/(q-1)$ into $\lambda=(p-1)(k-1)+1$ and exercising some algebra gives

\begin{equation*}%\label{eq:radialexponent2}
\lambda\left(\nu,\frac{q}{q-1}\right)=\frac{(\nu-1)\sqrt{(1-2\nu)(q-2)^2+\nu^2q^2}+(2-q)(1-2\nu)+\nu^2q}{2(q-1)(2\nu-1)}=k(\nu,q).
\end{equation*}

\bigskip
\noindent
Therefore, the same expression for the exponent continues to hold also for $q \in (1,2)$. Next, we consider the function $f_{\nu,p}$, defined in Case 2, in the extended domain $\phi \in [\frac{-\pi}{2\nu}, \frac{\pi}{\nu}]$.
Note that such extension is immediate (interpreting $\lim_{\underset{\theta\in(-\pi,\pi)}{\theta\to\pm\pi}}\phi(\theta)=\pm\frac{\pi}{\nu}$ in expression~\eqref{eq:implicitangel}) and that
%by~\cite{A86}, Case X on p. xxx]\komJ{inte hittat än},
the function
$(r, \phi) \to r^k f_{\nu,p}(\phi)$ is still $p$-harmonic in the extended domain.
Now, let $\bar u_{\nu,p}(x_1,x_2) = r^k f_{\nu,p}(\phi)$  for $r > 0$ and $\phi \in [0, \frac{\pi}{\nu}]$.
As we shall see, the stream function of  $\bar u_{\nu,p}$ is, up to rotation, the desired function.
In particular,  put $f_{\nu,p}(\phi)$ from \eqref{eq:losningv}  and $f'_{\nu,p}(\phi)$ from \eqref{eq:bounded derivative} in Lemma~\ref{le:konjugatlemma}, using $p=q/(q-1)$ and some algebra we arrive at
\begin{multline*}%\label{eq:stream3}
    g_{\nu,q}(\phi)=\frac{k\left(\nu,p\right)^{p-1}}{\lambda\left(\nu,\frac{q}{q-1}\right)}\left(1-\frac{\cos^2(\theta)}{ak(\nu,p)}\right)^{\frac{\left(k(\nu,p)-1\right)\left(p-1\right)}{2}}\sin(\theta)\\
    =\frac{k(\nu,p)^{p-1}}{k(\nu,q)}\left(1-\frac{\cos^2(\theta)}{\frac{1}{2-q}\left(\lambda(\nu,\frac{q}{q-1}) (q-1) + 2 - q\right)}\right)^{\frac{\lambda\left(\nu,\frac{q}{q-1}\right)-1}{2}}\sin(\theta)\\\
    = \frac{k(\nu,\frac{q}{q-1})^{\frac{1}{q-1}}}{k(\nu,q)}\left(1-\frac{\cos^2(\theta)}{\frac{q-1}{2-q}k(\nu,q) +1} \right)^{\frac{k(\nu,q)-1}{2}}\sin(\theta).
\end{multline*}

\bigskip

\noindent
Similarly, from the proof of Lemma~\ref{konjugatlemma} in Appendix~\ref{taste} we get
\begin{multline*}
 g_{\nu,q}'(\phi)=k(\nu,p) f_{\nu,p}(\phi)\left(k(\nu,p)^2f_{\nu,p}(\phi)^2+f'_{\nu,p}(\phi)^2\right)^{\frac{p-2}{2}}\\
 =k\left(\nu,\frac{q}{q-1}\right)^{\frac{1}{q-1}}\left(1-\frac{\cos^2(\theta)}{\frac{q-1}{2-q}k(\nu,q) +1} \right)^{\frac{k(\nu,q)-1}{2}}\cos(\theta).%=\\
% \gamma(v,q)k(v,q)\left(1-\frac{(2-q)\cos^2(\theta)}{k\left(v,\frac{q}{q-1}\right)\right)}\right)^{\frac{k(v,q)-1}{2}}\cos(\theta).
\end{multline*}

\bigskip

%As the $p$-Laplace equation is invariant under normalizations, we may ignore the factor in front and choose
%
%\begin{align*}
%g_{v,q}(\phi)=\left(1-\frac{\cos^2(\theta)}{\frac{q-1}{2-q}k(v,q) +1} \right)^{\frac{k(v,q)-1}{2}}\sin(\theta).
%\end{align*}
%
\noindent
Now since $\lim_{\nu\to 1/2}k(\nu,q)=\frac{q-1}{q}\leq k\left(\nu,q\right)$ and $0 < \frac{(q-1)^2}{(2-q)q}$ approaches zero only when $q \to 1$, we conclude
\begin{align}\label{eq:cont-deriavtion}
0< \kappa < 1 - \frac{\cos^2(\theta)}{\frac{q-1}{2-q}k(\nu,q) +1} < 1,
\end{align}

\bigskip
\noindent
whenever $q \in (1,2)$, where $\kappa$ depends only on $p$, and can be taken to be increasing in $p$. %which in fact can be taken decreasing in $p$ and independent of $\nu$.
Therefore, $|g_{\nu,q}| \leq c(\nu,p)$, $|g'_{\nu,q}| \leq c(\nu,p)$.

Now let $\tilde{f}_{\nu,p}(\phi)=g_{\nu,q}\left(\phi+\frac{\pi}{2\nu}\right)$. Then, since the $p$-Laplace equation is invariant under rotations (Lemma~\ref{le:scaling}), $u(x_1,x_2)=r^{k(\nu,q)}\tilde{f}_{\nu,p}(\phi)$ is a positive $q$-harmonic function in $\mathcal{S}_\nu$ satisfying the desired boundary conditions.
Moreover, $\tilde f'_{\nu,q}(\phi)$ is continuous by~\eqref{eq:cont-deriavtion} and only zero for $\phi = 0$ when restricting to $\mathcal{S}_\nu$. Therefore, we conclude that $|\tilde{f}'_{\nu,p}(\phi)|\geq c^{-1}$ whenever $\phi\in[-\frac{\pi}{2\nu},\frac{\pi}{2\nu}]\setminus[-\frac{\pi}{4\nu},\frac{\pi}{4\nu}]$. It also follows that $\tilde{f}_{\nu,p}(\phi) \geq c^{-1}$ whenever $\phi\in[-\frac{\pi}{4\nu},\frac{\pi}{4\nu}]$.
This completes the proof for the case $1<p<2$ and hence also the proof of Lemma~\ref{derivatalemma}.\end{proof}

 %\begin{comment}
 %Now for $0<\varepsilon<1$, we have for all $q\in(1+\varepsilon,2)$ that $\gamma(v,1+\varepsilon)\leq \frac{k\left(v,\frac{1+\varepsilon}{\varepsilon}\right)^{\frac{1}{\varepsilon}}}{k\left(\frac{1}{2},1+\varepsilon\right)}}=k\left(v,\frac{1+\varepsilon}{\varepsilon}\right)^{\frac{1}{\varepsilon}}(1+\varepsilon)=C(v,\varepsilon)$, hence $|g(\phi)|$ is less than or equal to a constant $C$ depending only on $v$ and $\varepsilon$. Similarly $|g'(\phi)|\leq \gamma(v,1+\varepsilon)k(v,1+\varepsilon)=D(v,\varepsilon)$.
 %\end{comment}

%%%%%%%%%%%%%%%%%%%%%%%%%%%%%%%%%%%%%%%%%%%%%%%%%%%%%%%%
%%%%%%%%%%%%%%%%%%%%%%%%%%%%%%%%%%%%%%%%%%%%%%%%%%%%%%%%
%%%%%%%%%%%%%%%%%%%%%%%%%%%%%%%%%%%%%%%%%%%%%%%%%%%%%%%%
%%%%%%%%%%%%%%%%%%%%%%%%%%%%%%%%%%%%%%%%%%%%%%%%%%%%%%%%
%%%%%%%%%%%%%%%%%%%%%%%%%%%%%%%%%%%%%%%%%%%%%%%%%%%%%%%%
%%%%%%%%%%%%%%%%%%%%%%%%%%%%%%%%%%%%%%%%%%%%%%%%%%%%%%%%
%%%%%%%%%%%%%%%%%%%%%%%%%%%%%%%%%%%%%%%%%%%%%%%%%%%%%%%%

\section{Estimates for $p$-harmonic measure}
\label{sec:est-p-meas}

\setcounter{theorem}{0}
\setcounter{equation}{0}

We will now state our growth estimate for $p$-harmonic measure in planar sectors. We postpone the proof to the end of the section.

\begin{theorem}\label{th:p-harmonic-measure}
Suppose that  $p \in (1,\infty]$, $\nu\in [1/2, \infty)$, $R>0$
and $\mathcal{S}_\nu \subset \mathbb{R}^2$ is the sector defined in \eqref{eq:def-S_v}.
Let $k(\nu,p)$ be the exponent in~\eqref{eq:radialexponent} and
let $\omega_p(x)$ be the $p$-harmonic measure of $\partial B(0,R) \cap \mathcal{S}_\nu$ at $x$
with respect to $B(0, R)\cap \mathcal{S}_\nu$.
Then there exists $c = c(\nu,p)$ such that
\begin{align*}
c^{-1}\,\left(\frac{|x|}{R}\right)^{k(\nu,p)}\,  \leq \omega_p(x)\, \leq c\, \left(\frac{|x|}{R}\right)^{k(\nu,p)},%,\quad e_1=(1,0).
\end{align*}
whenever $x \in B(0, R)\cap \mathcal{S}_{2\nu}$.
\end{theorem}

\noindent
%\begin{remark}
We remark that the upper bound of Theorem \ref{th:p-harmonic-measure} holds also in $B(0,R) \cap \mathcal{S}_{\nu}$ which follows from a comment in the end of the proof.
By Harnack's inequality in Lemma \ref{harnack} the lower bound also holds for any $x \in B(0, R)\cap \mathcal{S}_\nu$ but then the constant depends on the distance from $x$ to $\partial \mathcal{S}_\nu$.
Furthermore, by carefully tracing constants in the proof it can be shown that the final constant $c(\nu,p)$ in Theorem \ref{th:p-harmonic-measure} can be chosen independent of $p$ if $p$ is large, and the case $p = \infty$ in Theorem \ref{th:p-harmonic-measure} can be derived by taking the limit of the estimates for finite $p$.

For a final remark,
let $\bar\omega_p$ be the $p$-harmonic measure of
 $\partial B(0,R) \cap \mathcal{S}_{2\nu}$ at $x$
with respect to $B(0, R)\cap \mathcal{S}_\nu$.
Then there exists $c = c(\nu,p)$ such that
\begin{align}\label{eq:extra-theorem}
c^{-1}\,\left(\frac{|x|}{R}\right)^{k(\nu,p)}\,  \leq \bar\omega_p(x)\, \leq c\, \left(\frac{|x|}{R}\right)^{k(\nu,p)},%,\quad e_1=(1,0).
\end{align}

\bigskip
\noindent
whenever $x \in B(0, R/2)\cap \mathcal{S}_{2\nu}$.
To see that \eqref{eq:extra-theorem} holds we first observe that the upper bound is immediate from the comparison principle and Theorem \ref{th:p-harmonic-measure}.
To prove the lower bound we apply the boundary Harnack inequality in Lemma \ref{le:bhi-C11} to $\omega_p$ and $\bar\omega_p$ near the points of intersection
of $\partial B(0,R/2)$ and $\mathcal{S}_\nu$ and conclude that both functions vanishes at the same rate in a neighbourhood of these points. This, together with an application of H\"older continuity up to the boundary (Lemma \ref{le:holder}) near the point $(r,\phi) = (R,0)$ and Harnack's inequality toward the points of intersections, applied to $\bar \omega_p$, ensures $c(\nu,p) \bar\omega_p \geq \omega_p$ on $\partial B(0,R/2) \cap \mathcal{S}_\nu$. We now apply the comparison principle in
$B(0,R/2) \cap \mathcal{S}_\nu$ and Theorem \ref{th:p-harmonic-measure} to obtain the lower bound.
%\end{remark}

\bigskip

\noindent
{\bf Proof of Theorem \ref{th:p-harmonic-measure}.}
%\begin{proof}
Let $u(x) = u(r,\phi) = r^k f_{\nu,p}(\phi)$ be the $p$-harmonic function from Lemma \ref{derivatalemma} and
observe that $R^{-k}u(x)$ is $p$-subharmonic in $\mathcal{S}_\nu$ with boundary values $R^{-k}u(x) = 0$ on $\partial \mathcal{S}_\nu$ and $R^{-k}u(x) \leq 1$ on $\partial B(0,R) \cap  \mathcal{S}_\nu$, see $ii)$ in Lemma \ref{derivatalemma}.
This together with Definition \ref{def:p-hmeas} of the $p$-harmonic measure ensure that
\begin{equation*}
    \limsup_{x\rightarrow y}{R^{-k}u(x)}\leq\liminf_{x\rightarrow y}{\omega_p(x)},\quad \forall %x,y\in\mathbb{R}^2,\quad
    y\in \partial (B(0,R) \cap  \mathcal{S}_\nu),
\end{equation*}
and hence by the comparison principle (Lemma~\ref{jamforelseprin}) we obtain
\begin{align*}
R^{-k}u(x) \leq \omega_p(x), \quad \forall\,x\in B(0,R)\cap \mathcal{S}_\nu.
\end{align*}
Lemma \ref{derivatalemma} $iii)$ now implies
\begin{align}\label{lowerbound}
\omega_p(x)\geq R^{-k}u(x)=R^{-k}r^kf_{\nu,p}(\phi)=\left(\frac{|x|}{R}\right)^k f_{\nu,p}(\phi) \geq\frac{1}{c}\left(\frac{|x|}{R}\right)^k,
\end{align}
whenever $x \in B(0,R)\cap \mathcal{S}_{2\nu}$,
which establishes the lower bound.

\begin{figure}[!b]
\centering
\includegraphics[width=9cm,height=14cm,viewport=0 0 800 1300,clip]{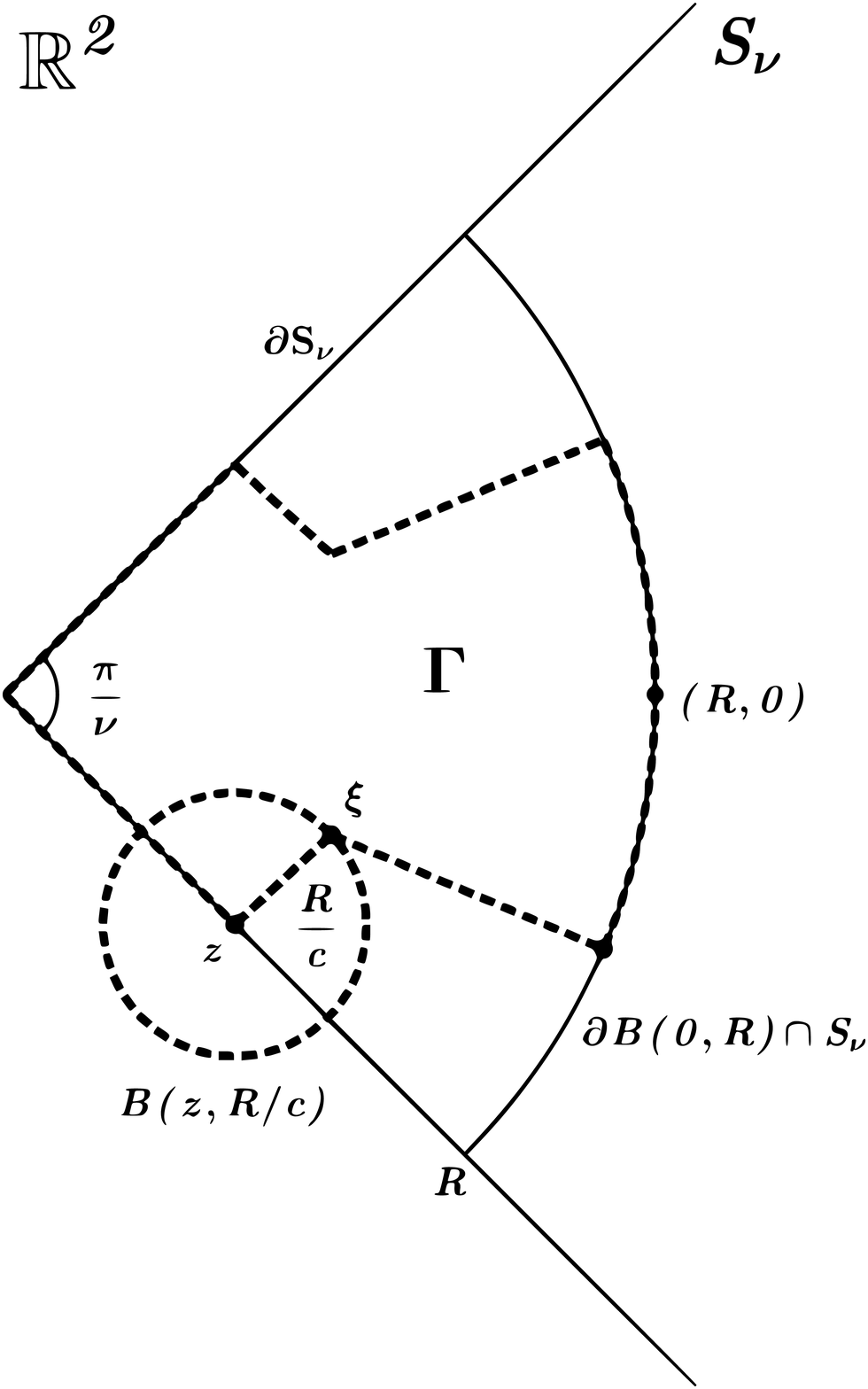}%{Fetgammafig2.eps}
\caption{Geometry in the proof of Theorem \ref{th:p-harmonic-measure}.}
\label{Gammafigur}
\end{figure}

We will now prove the upper bound.
As $\omega_p \approx 1$ at points near the intersections of $\partial B(0,R)$ and $\partial \mathcal{S}_\nu$ where $u \approx 0$ we chose to make comparison on a smaller domain not including these points.
Let $z \in \partial \mathcal{S}_\nu$ be the midpoint on the line from the origin to $(R, -\frac{\pi}{2\nu})$.
By the upper boundary growth estimate for $p$-harmonic functions in Lemma \ref{le:bhi-C11} there exists a constant $c = c(p)$ such that
\begin{align}\label{eq:upper-linear-growth}
\omega_p(x) \leq c\, \frac{d(x,\partial \mathcal{S}_\nu)}{R}\,\omega_p\left(a_R(w)\right)\leq c\, \frac{d(x,\partial \mathcal{S}_\nu)}{R} \quad \text{whenever} \quad x \in B(z,R/c) \cap \mathcal{S}_\nu.
\end{align}
Using the fact that $f_{\nu,p}'(\phi)$ does not vanish near $\partial \mathcal{S}_\nu$  (Lemma \ref{derivatalemma} $iii)$),
we also see that, for $c = c(\nu,p)$,
\begin{align}\label{eq-f-dont-vanish}
u(x)R^{-k}\geq\frac{1}{c} \frac{d(x,\partial \mathcal{S}_\nu)}{R} \quad \text{whenever} \quad x \in B(z,R/c) \cap \mathcal{S}_\nu.
\end{align}
Inequalities~\eqref{eq:upper-linear-growth} and \eqref{eq-f-dont-vanish} implies
\begin{align}\label{boll}
\omega_p(x)\leq c R^{-k} u(x) \quad \text{whenever} \quad x \in B(z,R/c) \cap \mathcal{S}_\nu.
\end{align}
Next, we also see from \eqref{eq-f-dont-vanish} and continuity of $u$ that,  for $c = c(\nu,p)$,
\begin{align}\label{eq:claim}
c R^{-k} u(\xi) \geq 1,
\end{align}
where $\xi = (R_\xi,\phi_\xi)$ is the point on the boundary of $B(z,R/c)$ with $d(\xi,\partial S_\nu) = R/c$.

Since $u$ and $\omega_p$ are continuous in $\mathcal{S}_\nu \cap B(0,R)$ and $0 < \omega_p < 1$ in $\mathcal{S}_\nu \cap B(0,R)$ we conclude
from \eqref{boll} and \eqref{eq:claim}
that
\begin{align}\label{eq:boundary-upper-est}
  \limsup_{x\to y}\,\omega_p(x)\leq\liminf_{x\to y} \, c R^{-k} u(x),\quad\forall\,  x\in \mathcal{S}_\nu \cap B(0,R) ,\,y \in\partial \Gamma.
\end{align}
Here, $\Gamma \subset \mathcal{S}_\nu$ is the open set bounded by the curve starting at the origin and reaching $z$ in the $r$-direction, then proceeding along a straight line to $\xi$ and from there to $\partial B(0,R)$ in $r$-direction, and proceeding in $\phi$-direction to the point $(R,0)$.
The rest of the curve, back to the origin, is the mirror of the above curve in the line $\phi = 0$ (see Figure~\ref{Gammafigur}).
Using \eqref{eq:boundary-upper-est} and the comparison principle (Lemma \ref{jamforelseprin}) we conclude that
\begin{align*}
\omega_p(x)\leq c R^{-k}u(x),\quad \forall x\in\Gamma.
\end{align*}
Since $B(0,R) \cap \mathcal{S}_{2\nu} \subset \Gamma$, at least if $|\phi_{\xi}|\geq \frac{\pi}{4\nu}$ which we may assume,  it follows that
\begin{align*}
\omega_p(x)\leq  c R^{-k}u(x) = c R^{-k}r^kf_{\nu,p}(\phi)\leq c \left(\frac{|x|}{R}\right)^{k},
\end{align*}
which establishes the upper bound in $B(0,R) \cap \mathcal{S}_{2\nu}$.
Observing that $B(0,R/2) \cap \mathcal{S}_{\nu} \subset \Gamma$ and that $\omega_p \leq 1$ we can conclude that the upper bound holds also whenever $x \in B(0,R) \cap \mathcal{S}_{\nu}$.
This together with the lower bound in \eqref{lowerbound} completes the proof of Theorem \ref{th:p-harmonic-measure}.
%\end{proof}
$\blacksquare$

%\begin{figure}[!hbt]\label{gammabild}
%\begin{center}	
%\includegraphics[width=8cm,height=8.5cm,viewport=220 35 460 330,clip]{gamma.eps}
%\end{center}
%\caption{The geometry of the proof of Theorem \ref{th:p-harmonic-measure}}
%\label{fig:geometry}
%\end{figure}

\bigskip

%\end{proof}
%$\blacksquare$\\

%\begin{remark}\label{rem:infty}
We remark that
in the case when $\nu \in [1/2, 1]$ and $p = \infty$,
giving $k(\nu,\infty) = 1$,
we can prove Theorem \ref{th:p-harmonic-measure} using the fact that infinity harmonic functions obey the comparison with cones principle, see \cite{CEG01}.
Indeed, the lower bound can be proved by making comparison with a cone function placed inside of $\mathcal{S}_\nu$, so that the circular base of the cone touches the boundary of $\mathcal{S}_\nu$ at the origin (always possible since $\nu \leq 1$).
The upper bound follows by comparison with a cone function placed so that its tip is at the origin.
This argument for proving an upper bound works for any $\nu \in [1/2, \infty)$ but it is optimal only when $\nu \leq 1$.

We also remark that for the cases when $p \in (1,\infty)$ and $\nu\in(\frac{1}{2},\infty)$ we may prove Theorem \ref{th:p-harmonic-measure} by scaling and an application of the boundary Harnack inequality in Lemma \ref{le:lewis-nyström} (given in Section \ref{sec:unique} below) by taking $v(x)$ as the $p$-harmonic function in Lemma \ref{derivatalemma} and $u(x)$ as the $p$-harmonic measure.
This works because the sector $\mathcal{S}_\nu$ is a Lipschitz domain as long as $\nu \neq \frac{1}{2}$.
%\end{remark}

%%%%%%%%%%%%%%%%%%%%%%%%%%%%%%%%%%%%%%%%%%%%%%%%%%%%%%%%
%%%%%%%%%%%%%%%%%%%%%%%%%%%%%%%%%%%%%%%%%%%%%%%%%%%%%%%%
%%%%%%%%%%%%%%%%%%%%%%%%%%%%%%%%%%%%%%%%%%%%%%%%%%%%%%%%
%%%%%%%%%%%%%%%%%%%%%%%%%%%%%%%%%%%%%%%%%%%%%%%%%%%%%%%%
%%%%%%%%%%%%%%%%%%%%%%%%%%%%%%%%%%%%%%%%%%%%%%%%%%%%%%%%
%%%%%%%%%%%%%%%%%%%%%%%%%%%%%%%%%%%%%%%%%%%%%%%%%%%%%%%%
%%%%%%%%%%%%%%%%%%%%%%%%%%%%%%%%%%%%%%%%%%%%%%%%%%%%%%%%

\section{Estimates for $p$-sub- and $p$-superharmonic functions}
\label{sec:est-p-harmonic}

\setcounter{theorem}{0}
\setcounter{equation}{0}

In this section we state and prove some Corollaries of Theorem \ref{th:p-harmonic-measure} giving estimates of $p$-sub- and $p$-superharmonic functions in domains related to planar sectors.

\begin{corollary}\label{cor:subsuper}
Suppose that $p \in (1,\infty]$, $\nu \in [1/2,\infty)$, $R > 0$,
 $\Omega \in \mathbb{R}^2$ is a domain, $w \in \partial \Omega$ and that
$\Omega \cap B(w,R)$ is contained in a planar sector with apex
$w$ and aperture angle $\frac{\pi}{\nu}$.
Let $u$ be a $p$-subharmonic function in $\Omega$,
satisfying $u \leq 0$ on $\partial \Omega \cap B(w,R)$.  Then
\begin{align}\label{eq:first-cor}
%\begin{itemize}
    %\item
u(x) \leq c M\left(\frac{|x-w|}{R}\right)^{k(\nu,p)}
\end{align}

\noindent
whenever $x \in \Omega\cap B(w, R)$,
$M = \max\left\{0, \sup_{\partial B(w,R)\cap \Omega} u\right\}$
and where $c=c(\nu,p)$ is the constant in Theorem \ref{th:p-harmonic-measure}.

Suppose now instead that $\mathcal{C}_\nu \cap B(w,R)$ is contained in $\Omega$, where $\mathcal{C}_\nu$ is a cone with apex
$w$ and aperture angle $\frac{\pi}{\nu}$,
and that $v$ is a nonnegative $p$-superharmonic function in $\Omega$,
satisfying $v \geq 0$ on $\partial \Omega \cap B(w,R)$.
Then
\begin{align}\label{eq:first-cor2}
%\qquad \text{and} \qquad
    %\item
 c^{-1} m \left(\frac{|x-w|}{R}\right)^{k(\nu,p)} \leq v(x)
%\end{itemize}
%c_1^{-1} u(x) \leq |x|^{k(v,p)} \leq c_2 v(x)
\end{align}
whenever $x \in B(w, R/2) \cap \mathcal{C}_{2v}$,
$m = \inf_{\partial B(w,R)\cap \mathcal{C}_{2v}} v$
and where $c=c(\nu,p)$ is the constant in inequality  \eqref{eq:extra-theorem}.%Theorem \ref{th:p-harmonic-measure}.
%whenever $x \in \mathcal{S}_v \cap B(0, r/c)$.
%Here, $c_1$ may depend on $p$, $v$ and $\sup_{x\in B(0,R)\cap \mathcal{S}_v} u$,
%and  $c_2$ may depend on $p$, $v$ and $\inf_{x\in B(0,R)\cap \mathcal{S}_v} v$
\end{corollary}

\noindent
%{\bf Proof of Corollary \ref{cor:subsuper}.}
\begin{proof}
%\noindent
%\begin{proof}
We begin by proving \eqref{eq:first-cor}.
Thanks to Lemma \ref{le:scaling} we change coordinates so that $w = (0,0)$ and the domain $\Omega$ is contained in $\mathcal{S}_\nu$.
Let $\omega_p$ be the $p$-harmonic measure in Theorem \ref{th:p-harmonic-measure}.
We can conclude that
$$
M \omega_p(x) \leq c M \left(\frac{|x|}{R}\right)^{k(v,p)}
$$
whenever $x\in B(0,R)\cap \Omega$ where
$$
M = \max\left\{0, \sup_{\partial B(0,R)\cap \Omega} u\right\}.
$$
Moreover,
$M \omega_p$ is $p$-harmonic with boundary values dominating $u$ on $\partial (B(0,R)\cap \Omega)$.
Therefore, by the comparison principle in Lemma \ref{jamforelseprin} we have
$$
u(x) \leq M \omega_p(x)
$$
whenever $x\in B(0,R)\cap \Omega$ and
the first inequality in Corollary \ref{cor:subsuper} follows by returning to the original coordinates.

To prove \eqref{eq:first-cor2},
change coordinates so that $w = (0,0)$, $\mathcal{C}_\nu = \mathcal{S}_\nu$, and
let $\bar \omega_p$ be the $p$-harmonic measure in \eqref{eq:extra-theorem}. %Theorem \ref{th:p-harmonic-measure}.
We can conclude that
$$
c^{-1} m \left(\frac{|x|}{R}\right)^{k(\nu,p)} \leq m \,\bar \omega_p(x)
$$
whenever $x\in B(0,R/2)\cap \mathcal{S}_{2\nu}$ where
$
m = \inf_{\partial B(0,R)\cap \mathcal{S}_{2\nu}} v.
$
It follows that $m\, \bar \omega_p$ is $p$-harmonic with boundary values dominated by $v$ on $\partial (B(0,R)\cap \mathcal{S}_\nu)$.
Therefore, by the comparison principle in Lemma \ref{jamforelseprin} we have
$$
m\, \bar\omega_p(x) \leq v(x)
$$
whenever $x\in B(0,R)\cap \mathcal{S}_\nu$ and
the second inequality in Corollary \ref{cor:subsuper} follows by returning to the original coordinates.\end{proof}

\bigskip

Corollary \ref{cor:subsuper} implies growth estimates for $p$-sub- and $p$-superharmonic functions near the boundary of a large class of planar domains.
Consider e.g. a domain $\Omega \subset \mathbb{R}^2$ having a sharp outwardly pointed cusp with apex $w$.
Then, in a neighborhood of $w$, the domain will be contained in a planar sector with small aperture angle and apex at $w$, and as the neighborhood shrinks the aperture angle of the sector also shrinks, i.e. $\nu \to \infty$ in Corollary \ref{cor:subsuper}. Since $k(\nu,p) \to \infty$ as $\nu\to \infty$ it follows from \eqref{eq:first-cor}
that if a $p$-subharmonic function $u(x)$ takes nonpositive boundary values in the neighborhood of the apex $w$, then the rate of convergence to zero, as $x$ approaches the apex, is faster than any power of $|x-w|$.
Indeed, for any $N > 0$ it holds that
%
%\begin{align}
%\lim_{\Omega \ni x\to w} u(x) |x-w|^N = 0,
%\end{align}
%

\begin{align*}
\underset{x\in\Omega}{\underset{x\to w}{\limsup}}\;\frac{u(x)}{|x-w|^N}\leq 0,
\end{align*}

\bigskip
\noindent
which is a result proved already in \cite[Theorem 3]{rysk}.
Using \eqref{eq:first-cor2} in Corollary \ref{cor:subsuper} we now derive a similar estimate for $p$-superharmonic functions in domains $\Omega \subset \mathbb{R}^2$ having a sharp inwardly pointed cusp at $w$.
Indeed, in such case the domain will, in a neighborhood of $w$, contain a planar sector with large aperture angle ($\nu$ close to $\frac{1}{2}$ in Corollary \ref{cor:subsuper}) and as the neighborhood shrinks we may let $\nu \to \frac{1}{2}$ implying $k(\nu,p) \to \frac{p-1}{p}$.
It follows from \eqref{eq:first-cor2}
that if a positive $p$-superharmonic function $ v(x)$ takes nonnegative boundary values in the neighborhood of the apex $w$, then the rate of convergence to zero, as $x$ approaches the apex, is slower than $|x-w|^N$ whenever $N > \frac{p-1}{p}$.
Going to the limit with $N$ implies
%
%\begin{align}
%\lim_{\Gamma \ni x\to w} \frac{u(x)}{|x-w|^N} = 0, \qquad \text{where}\qquad \Gamma = \{x \in \Omega : d(x,\partial \Omega) \geq |x-w|\}.
%\end{align}
%
\begin{align*}
\underset{x\in\Lambda}{\underset{x\to w}{\liminf}}\;\frac{ v(x)}{|x-w|^\frac{p-1}{p}} > 0, \quad \text{where}\quad \Lambda = \{x \in \Omega : d(x,\partial \Omega) \geq |x-w|\}.
\end{align*}

\bigskip

%%%% Letter from Mikhail Surnachev about russian paper \cite{rysk}

%The author studies the elliptic p-Laplacian with the homogeneous boundary conditions in the cone $K(l)=\{x: 0\leq\theta\leq l\}$ where $\cos \theta=x_n/|x|$. He constructs on this cone solutions of the form $u(x)=|x|^\lambda f(\theta)$ and proves that $\lambda$ as the function of the opening angle of the cone has the asymptotics $\lambda(l)=\pm L l^{-1}+O(1)$. This asymptotics is the result of Theorem 2 of the paper, whose statement is at the bottom of page 4 and continues to page 5. It says that $\lambda(l)$ has asymptotics (12) where the number $L$ is the first zero of the solution of the Cauchy problem (13) with the initial data (14). On page 7 the author writes that if $n=2$ one can write down an explicit quadratic equation for $\lambda$ (approx in the center of the page), which gives the value of $L=\pm \pi p/(4(p-1))$. In the end of the paper the author uses these results to analyze the boundary behavior of a p-harmonic function in the domain with a sharp outwardly pointed peak. Namely, he shows that if a solution takes zero boundary value in the neighborhood of the apex $x_0$ of a sharp peak than the rate of convergence to this boundary value is faster than any power of $|x-x_0|$.

Using Corollary \ref{cor:subsuper} we can also derive
the boundary Harnack's inequality for positive $p$-harmonic functions vanishing on a portion of the boundary of a planar sector:
Suppose that $p \in (1,\infty]$, $\nu \in [1/2,\infty)$, $R > 0$ and $w \in \partial \mathcal{S}_\nu$.
Suppose also that $u_1$ and $u_2$ are positive $p$-harmonic functions in $\mathcal{S}_\nu \cap B(w, 2R)$, satisfying $u_1 = 0 = u_2$ on $\partial \mathcal{S}_\nu \cap B(w, 2R)$.
Then, if $w$ is the apex of the sector there exists $c = c(\nu,p)$ such that, for $i = 1,2$,
\begin{align}\label{eq:second-cor}
 c^{-1} \left(\frac{|x|}{R}\right)^{k(\nu,p)} \leq \frac{u_i(x)}{u_i(R,0)} \leq c  \left(\frac{|x|}{R}\right)^{k(\nu,p)} \quad \text{and} \quad c^{-1} \frac{u_1(R,0)}{u_2(R,0)} \leq \frac{u_1(x)}{u_2(x)} \leq c   \frac{u_1(R,0)}{u_2(R,0)} ,
\end{align}

\bigskip
\noindent
whenever $x \in \mathcal{S}_{2\nu} \cap B(0, R/2)$, and where we have used the polar coordinates notation $u_i(x) = u_i(r,\phi)$.
To derive \eqref{eq:second-cor} from Corollary \ref{cor:subsuper} we observe that Lemma \ref{le:bhi-C11} and Harnack's inequality  (Lemma \ref{harnack}) implies the existence of a constant
$c = c(\nu,p)$ such that, for $i = 1,2$,
$$
\frac{1}{c}\, u_i(R,0) \leq \inf_{\partial B(0,R)\cap \mathcal{S}_{2\nu}} u_i \qquad \text{and} \qquad \sup_{\partial B(0,R)\cap \mathcal{S}_\nu} u_i \leq c\, u_i(R,0).
$$

\bigskip

%Alternatively, for $p \in (1,\infty)$, \eqref{eq:second-cor} follows by an application of Lemma \ref{le:lewis-nyström} putting $v(x)$ to the $p$-harmonic function in Lemma \ref{eq:radialexponent}.

\noindent
The left inequality in \eqref{eq:second-cor} states that any positive $p$-harmonic function, vanishing on the boundary of the sector,
must vanish at the same rate as the distance to the apex to the power of $k(\nu,p)$.
The right inequality in Lemma \ref{le:bhi-C11} - which is an immediate consequence of the left inequality - is usually referred to as a boundary Harnack inequality and states that any two $p$-harmonic functions must vanish at the same rate.
If $w$ is not the apex of the sector $\mathcal{S}_\nu$
then in a neighbourhood of $w$ the boundary is a line and the estimates in \eqref{eq:second-cor} are well known to hold with $k(\nu,p) = 1$.
In particular, such result is given in Lemma \ref{le:bhi-C11} and was proved in  %Aikawa--Kilpel\"ainen--Shanmugalingam--Zhong~
\cite{aikawa} for $C^{1,1}$-domains and  in \cite{LN07,LN08,LN10} for Lipschitz and Reifenberg flat doamins.

%%%%%%%%%%%%%%%%%%%%%%%%%%%%%%%%%%%%%%%%%%%%%%%%%%%%%%%%
%%%%%%%%%%%%%%%%%%%%%%%%%%%%%%%%%%%%%%%%%%%%%%%%%%%%%%%%
%%%%%%%%%%%%%%%%%%%%%%%%%%%%%%%%%%%%%%%%%%%%%%%%%%%%%%%%
%%%%%%%%%%%%%%%%%%%%%%%%%%%%%%%%%%%%%%%%%%%%%%%%%%%%%%%%
%%%%%%%%%%%%%%%%%%%%%%%%%%%%%%%%%%%%%%%%%%%%%%%%%%%%%%%%
%%%%%%%%%%%%%%%%%%%%%%%%%%%%%%%%%%%%%%%%%%%%%%%%%%%%%%%%
%%%%%%%%%%%%%%%%%%%%%%%%%%%%%%%%%%%%%%%%%%%%%%%%%%%%%%%%

\section{A sharp Phragmen-Lindelöf theorem}
\label{sec:phragmen}

\setcounter{theorem}{0}
\setcounter{equation}{0}

In this section we will prove sharp lower growth estimates of
$p$-subharmonic functions in planar sectors.
To state our theorem, let $\Omega \subset \mathbb{R}^2$ be a domain contained in a planar sector.
Assume without loss of generality (thanks to Lemma \ref{le:scaling}) that the sector is $\mathcal{S}_\nu$ given in \eqref{eq:def-S_v} and define
%
%\begin{align*}
%M(R)= \sup_{x\in \mathcal{S}_v,\,|x|=R} u(x).
%\end{align*}

\begin{align*}
M(R)=\sup_{\partial B(0,R) \cap \Omega} u,
\end{align*}

%\begin{align*}
%M(R)=\underset{x\in S_v}{\underset{|x|=R}{\sup}}\; u(x).
%\end{align*}

\noindent
for $R > 0$.
Using the estimates of $p$-harmonic measure in Theorem \ref{th:p-harmonic-measure} we obtain the following version of the Phragmen-Lindel\"of theorem:
\begin{theorem}\label{th:phragmen}
Suppose that $p \in (1, \infty]$, $\nu \in [1/2, \infty)$
and that $u$ is a $p$-subharmonic function in a domain $\Omega \subset \mathcal{S}_\nu$ satisfying
\begin{align*}
\limsup_{x \to y} u(x) \leq 0 \quad \text{for each}\; y \in \partial \Omega.
\end{align*}
Then either $u \leq 0$ in $\Omega$ or it holds that
\begin{align*}
\liminf_{R \to \infty} \frac{M(R)}{R^{k(\nu,p)}} > 0,
\end{align*}
where $k(\nu,p)$ is the exponent in \eqref{eq:radialexponent}.
%
%\begin{equation*}%\label{eq:radialexponent}
%k(v,p)=\frac{\sqrt{(1-2v)(p-2)^2+v^2p^2}(v-1)+(2-p)(1-2v)+v^2p}{2(p-1)(2v-1)}.
%\end{equation*}
%
\end{theorem}
In case $\Omega = \mathcal{S}_\nu$ then the $p$-harmonic function from Lemma \ref{derivatalemma} shows that the growth estimate in Theorem \ref{th:phragmen} is sharp.
The proof uses the following well known Phragmen-Lindel\"of principle which can be found in a more general form in~\cite[11.11]{HKM}, and is a key to the study of the
behaviour of $M(R)$.

\begin{lemma}\label{le:principle}
Let $p \in (1, \infty]$, $\nu \in [1/2, \infty)$, $u$ be as in Theorem \ref{th:phragmen}, and suppose
for each $R>0$ that $v(x)$ is $p$-superharmonic in  $\mathcal{S}_\nu$ with
\begin{align*}
\lim_{x \to y} v(x) = 1, \quad  y \in \partial B(0,R) \cap \mathcal{S}_\nu.
\end{align*}
Then either $u \leq 0$ in $\mathcal{S}_\nu$ or it holds that
\begin{align*}
\liminf_{R \to \infty} (M(R) v(x)) > 0,
\end{align*}
for any $x \in \mathcal{S}_\nu$.
\end{lemma}

\noindent
\begin{proof}
This follows from The Phragmen-Lindel\"of principle~\cite[11.11]{HKM}.\end{proof}

\bigskip

\noindent
{\bf Proof of Theorem \ref{th:phragmen}.}
The result follows from our estimates in Theorem \ref{th:p-harmonic-measure} by taking $v$ in Lemma \ref{le:principle} as the $p$-harmonic measure.%\end{proof}
$\blacksquare$

%%%%%%%%%%%%%%%%%%%%%%%%%%%%%%%%%%%%%%%%%%%%%%%%%%%%%%%%%%%%%%%%%%%%%%%%%%%%%%%%%%%%
%%%%%%%%%%%%%%%%%%%%%%%%%%%%%%%%%%%%%%%%%%%%%%%%%%%%%%%%%%%%%%%%%%%%%%%%%%%%%%%%%%%%
%%%%%%%%%%%%%%%%%%%%%%%%%%%%%%%%%%%%%%%%%%%%%%%%%%%%%%%%%%%%%%%%%%%%%%%%%%%%%%%%%%%%
%%%%%%%%%%%%%%%%%%%%%%%%%%%%%%%%%%%%%%%%%%%%%%%%%%%%%%%%%%%%%%%%%%%%%%%%%%%%%%%%%%%%
%%%%%%%%%%%%%%%%%%%%%%%%%%%%%%%%%%%%%%%%%%%%%%%%%%%%%%%%%%%%%%%%%%%%%%%%%%%%%%%%%%%%
%%%%%%%%%%%%%%%%%%%%%%%%%%%%%%%%%%%%%%%%%%%%%%%%%%%%%%%%%%%%%%%%%%%%%%%%%%%%%%%%%%%%
%%%%%%%%%%%%%%%%%%%%%%%%%%%%%%%%%%%%%%%%%%%%%%%%%%%%%%%%%%%%%%%%%%%%%%%%%%%%%%%%%%%%
%%%%%%%%%%%%%%%%%%%%%%%%%%%%%%%%%%%%%%%%%%%%%%%%%%%%%%%%%%%%%%%%%%%%%%%%%%%%%%%%%%%%
%%%%%%%%%%%%%%%%%%%%%%%%%%%%%%%%%%%%%%%%%%%%%%%%%%%%%%%%%%%%%%%%%%%%%%%%%%%%%%%%%%%%
%%%%%%%%%%%%%%%%%%%%%%%%%%%%%%%%%%%%%%%%%%%%%%%%%%%%%%%%%%%%%%%%%%%%%%%%%%%%%%%%%%%%

\section{Uniqueness of $p$-harmonic functions in sectors}
\label{sec:unique}

\setcounter{theorem}{0}
\setcounter{equation}{0}

It is well known that a positive $p$-harmonic function in the halfspace $\mathbb{R}^n_+$,
vanishing on the boundary, %(an $(n-1)$-dimensional plane)
must be a multiple of the distance to the boundary.
In case of $n = 2$,
the following theorem generalizes this result to planar sectors.

\begin{theorem}\label{th:unique}
Let $p \in (1,\infty)$, $\nu \in (1/2,\infty)$
and suppose that $u$ is a positive $p$-harmonic function in the sector $\mathcal{S}_\nu$.
Suppose also that $u = 0$ on $\partial \mathcal{S}_\nu$.
Then there exists a constant $c$ such that $u = c\, r^k f_{p,\nu}(\phi)$
where $k = k(\nu,p)$ is as in \eqref{eq:radialexponent} and $r^k f_{p,\nu}(\phi)$ is the $p$-harmonic function in Lemma \ref{derivatalemma}. \end{theorem}

\noindent
The proof uses the following boundary estimate from \cite{LN10},
valid for $p \in (1, \infty)$, stating that the ratio of two positive $p$-harmonic functions,
 both vanishing on a portion of a Lipschits boundary,
is H\"older continuous near the boundary.
The reason for excluding $\nu = 1/2$ is that then $\mathcal{S}_\nu$ fails to be Lipschitz.

\begin{lemma}\label{le:lewis-nyström}
Let $\Omega \subset \mathbb{R}^n$ be a bounded Lipschitz domain with constant $M$.
Given $p \in (1,\infty), w \in \partial\Omega$,
and $0 < r \leq r_0$ for some $r_0 < \infty$,
suppose that $u$ and $v$ are positive $p$-harmonic functions in $\Omega \cap B(w, r)$. Assume also that $u$ and $v$ are continuous
in $\overline\Omega \cap B(w, r)$ and that $u = 0 = v$ on $\Omega \cap \partial B(w, r)$.
Under these assumptions there exist $c \in (1,\infty)$ and $\alpha \in (0,1)$,
both depending only on $p, n$ and $M$, such that if
$y_1,y_2 \in \Omega \cap B(w, r/c)$, then
$$
\left| \log\frac{u(y_1)}{v(y_1)} - \log \frac{u(y_2)}{v(y_2)} \right| \leq c \left(\frac{|y_1-y_2|}{r}\right)^\alpha.
$$
\end{lemma}

\noindent
\begin{proof}
See \cite[Theorem 2]{LN10}.\end{proof}\\

\noindent
%\begin{proof}
{\bf Proof of Theorem \ref{th:unique}.}
%We intend to apply Lemma \ref{le:lewis-nyström}.
Let $\mathcal{S}_\nu$ and $u$ be as in the theorem and consider the bounded sector $\Omega = \mathcal{S}_\nu \cap \{r < 1\}$. It is clear that $\Omega$ is a bounded Lipschitz domain for any $\nu \in (\frac{1}{2},\infty)$.
Define the scaled function
$%\begin{align*}
u_1(x) = u(R x)
$. %\end{align*}
Then, since $u$ is $p$-harmonic in $B(0,R) \cap \mathcal{S}_\nu$ it follows by Lemma \ref{le:scaling} that $u_1$ is $p$-harmonic in $\Omega$.
Let $v_1$ be the explicit $p$-harmonic function in Lemma \ref{derivatalemma}, scaled in the same way as $u$.
Then $v_1$ is also $p$-harmonic in $\Omega$.
As $\Omega$ is a bounded Lipschitz domain with Lipschitz constant $M$ depending only on $\nu$, and since $u_1$ and $v_1$ are zero on the sides of the sector $\Omega$, we deduce from Lemma  \ref{le:lewis-nyström}, with $\omega = 0$, and $r = r_0 = 1$, that
$$
\left| \log\frac{u_1(y_1)}{v_1(y_1)} - \log \frac{u_1(y_2)}{v_1(y_2)} \right| \leq c \left(|y_1-y_2|\right)^\alpha,
$$
whenever $y_1,y_2 \in \Omega \cap B(0, 1/c)$ and $c = c(\nu,p)$.
Let $x^1,x^2$ be arbitrary points in $\mathcal{S}_\nu$.
Pick $R$ so large
that $x^1, x^2 \in S_\nu \cap B(0,R/c)$ where $c$ is from the above display. In the scaled domain, these points are $x_1 = x^1/R, x_2 = x^2/R$ and they end up in $\Omega \cap B(0,1/c)$. Thus
$$
\left| \log\frac{u_1(x_1)}{v_1(x_1)} - \log \frac{u_1(x_2)}{v_1(x_2)} \right| \leq c \left(|x_1-x_2|\right)^\alpha =
c \left(\frac{|x^1-x^2|}{R}\right)^\alpha.
$$
As $R$ can be taken arbitrary large we may send $R \to \infty$ and thereby deduce, since also $x^1$ and $x^2$ were arbitrary, that  ${u_1}/{v_1}$ must be constant and therefore $u_1 = c \,v_1$ for some constant $c$. Scaling back concludes the proof.
$\blacksquare$
%\end{proof}

'%When $v = 1/2$ the boundary is not Lipschitz at the origin and we therefore proceed along the following arguments:
%\end{proof}
%$\blacksquare$\\

%We remark that arguments in the proof of Theorem \ref{th:unique} is valid in $n \geg 2$ dimensions as well. The reason for   \komN{for $x \in \mathbb{R}^2$ but this argument works also in $n$-dim. Let us check some refs for explicit $p$-harmonic functions in $\mathbb{R}^n$ later on}

%%%%%%%%%%%%%%%%%%%%%%%%%%%%%%%%%%%%%%%%%%%%%%%%%%%%%%%%%%%%%%%%%%%%%%%%%%%%%%%%%%%%
%%%%%%%%%%%%%%%%%%%%%%%%%%%%%%%%%%%%%%%%%%%%%%%%%%%%%%%%%%%%%%%%%%%%%%%%%%%%%%%%%%%%
%%%%%%%%%%%%%%%%%%%%%%%%%%%%%%%%%%%%%%%%%%%%%%%%%%%%%%%%%%%%%%%%%%%%%%%%%%%%%%%%%%%%
%%%%%%%%%%%%%%%%%%%%%%%%%%%%%%%%%%%%%%%%%%%%%%%%%%%%%%%%%%%%%%%%%%%%%%%%%%%%%%%%%%%%
%%%%%%%%%%%%%%%%%%%%%%%%%%%%%%%%%%%%%%%%%%%%%%%%%%%%%%%%%%%%%%%%%%%%%%%%%%%%%%%%%%%%
%%%%%%%%%%%%%%%%%%%%%%%%%%%%%%%%%%%%%%%%%%%%%%%%%%%%%%%%%%%%%%%%%%%%%%%%%%%%%%%%%%%%
%%%%%%%%%%%%%%%%%%%%%%%%%%%%%%%%%%%%%%%%%%%%%%%%%%%%%%%%%%%%%%%%%%%%%%%%%%%%%%%%%%%%
%%%%%%%%%%%%%%%%%%%%%%%%%%%%%%%%%%%%%%%%%%%%%%%%%%%%%%%%%%%%%%%%%%%%%%%%%%%%%%%%%%%%
%%%%%%%%%%%%%%%%%%%%%%%%%%%%%%%%%%%%%%%%%%%%%%%%%%%%%%%%%%%%%%%%%%%%%%%%%%%%%%%%%%%%
%%%%%%%%%%%%%%%%%%%%%%%%%%%%%%%%%%%%%%%%%%%%%%%%%%%%%%%%%%%%%%%%%%%%%%%%%%%%%%%%%%%%

\section{Extension to $n$-dimensional cones when $p = \infty$}
\label{sec:inf}

Assume $n \geq 2$
%a domain $\Omega \subseteq \mathbb{R}^n$ is rotationally invariant around an axis $\ell$.
%Let $\Omega_{\ell}$ be the intersection between $\Omega$ and a two-dimensional plane
%containing $\ell$.
%We define the $n$-dimensional cone $\mathcal{S}_v^n$ as $\Omega$ when $\Omega_{\ell} \equiv \mathcal{S}_v^n$ where $\mathcal{S}_v$ is the sector in \eqref{eq:def-S_v}.
and define the $n$-dimensional cone $\mathcal{S}_\nu^n$ as a domain being rotationally invariant around the $x_1$ axis and of which its intersection with any two-dimensional plane
containing the $x_1$ axis equals $\mathcal{S}_\nu$ (modulo rotation).
%\begin{align}\label{eq:def-S_v-n}
%\mathcal{S}_v^n = \left\{ (x_1, x_2, \dots, x_n) \in  \mathbb{R}^n \setminus \{(0, 0)\}; |\phi| <
%\frac{\pi}{2v}\right\} \quad \text{where $v \geq \frac{1}{2}$}
%\end{align}
%
Recall that the infinity-Laplace equation is invariant under rotations, scaling and translations (Lemma \ref{le:scaling}) and hence the following corollary applies to any $n$-dimensional cone.
In the case $p = \infty$ we have the following extension of our Theorems from planar domains into $\mathbb{R}^n$:
\begin{corollary} \label{cor:inf-hmeas}
Suppose that $n \geq 2$, $\nu\in[1/2,\infty)$ and that $p = \infty$.
Then Theorem \ref{th:p-harmonic-measure}, Corollary \ref{cor:subsuper} and Theorem \ref{th:phragmen} generalize to the corresponding $n$-dimensional setting.
In particular, these results hold also when the two-dimensional cone $\mathcal{S}_\nu$ is replaced by the $n$-dimensional cone  $\mathcal{S}_\nu^n$, $\Omega \subset \mathbb{R}^n$, and $k(\nu,p) = k(\nu,\infty)$ is as in \eqref{eq:kinfhej}.
\end{corollary}

\noindent
\begin{proof}
Corollary \ref{cor:subsuper} and Theorem \ref{th:phragmen} follow from Theorem \ref{th:p-harmonic-measure} by standard arguments which are valid in $\mathbb{R}^n$ as well.
Therefore, we focus on the extension of Theorem \ref{th:p-harmonic-measure} from two to $n$-dimensions.

Suppose that $\omega = \omega_\infty$ satisfies the assumptions in the theorem but in $n$-dimensions, $n > 2$.
Then $\omega$ is $\infty$-harmonic in $B(0,R) \cap \mathcal{S}_\nu^n$. We will show that by symmetry, $\omega$ is also $\infty$-harmonic in the two-dimensional sector $B(0,R) \cap \mathcal{S}_\nu$ and therefore the result remains.
Assume that $\omega \in C^2(\Omega)$, otherwise, we switch to a $C^2$-function through the definition of viscosity solutions.
%Since the $\infty$-Laplacian is invariant under rotations and translations, we also assume that the axis $\ell$ in the coincides with the $x_1$-axis and that $\Omega_{\ell}$ is contained in the $x_1 x_2$-plane.
By symmetry of the bounded domain $B(0,R) \cap \mathcal{S}_\nu^n$, symmetry of the boundary conditions, and by the fact that the $\infty$-harmonic measure is unique, we conclude that  $\omega_{x_3} = \omega_{x_4} = \dots = \omega_{x_n} = 0$ on the two-dimensional cone $B(0,R) \cap \mathcal{S}_\nu$ and hence
\begin{equation*}
\Delta_{\infty} \omega  = \sum_{i, j=1}^{n} \omega_{x_i}\omega_{x_j}\omega_{x_ix_j} = \omega_{x_1}^2\omega_{x_1x_1} + 2\omega_{x_1}\omega_{x_2}\omega_{x_1x_2} + \omega_{x_2}^2\omega_{x_2x_2} = 0.
\end{equation*}
Thus, $\omega$ is $\infty$-harmonic in $B(0,R)\cap\mathcal{S}_\nu \subseteq \mathbb{R}^2$ and we conclude Corollary~\ref{cor:inf-hmeas}.
\end{proof}

%%%%%%%%%%%%%%%%%%%%%%%%%%%%%%%%%%%%%%%%%%%%%%%%%%%%%%%%%%%%%%%%%%%%%%%%%%%%%%%%%%%%
%%%%%%%%%%%%%%%%%%%%%%%%%%%%%%%%%%%%%%%%%%%%%%%%%%%%%%%%%%%%%%%%%%%%%%%%%%%%%%%%%%%%
%%%%%%%%%%%%%%%%%%%%%%%%%%%%%%%%%%%%%%%%%%%%%%%%%%%%%%%%%%%%%%%%%%%%%%%%%%%%%%%%%%%%
%%%%%%%%%%%%%%%%%%%%%%%%%%%%%%%%%%%%%%%%%%%%%%%%%%%%%%%%%%%%%%%%%%%%%%%%%%%%%%%%%%%%
%%%%%%%%%%%%%%%%%%%%%%%%%%%%%%%%%%%%%%%%%%%%%%%%%%%%%%%%%%%%%%%%%%%%%%%%%%%%%%%%%%%%
%%%%%%%%%%%%%%%%%%%%%%%%%%%%%%%%%%%%%%%%%%%%%%%%%%%%%%%%%%%%%%%%%%%%%%%%%%%%%%%%%%%%
%%%%%%%%%%%%%%%%%%%%%%%%%%%%%%%%%%%%%%%%%%%%%%%%%%%%%%%%%%%%%%%%%%%%%%%%%%%%%%%%%%%%
%%%%%%%%%%%%%%%%%%%%%%%%%%%%%%%%%%%%%%%%%%%%%%%%%%%%%%%%%%%%%%%%%%%%%%%%%%%%%%%%%%%%
%%%%%%%%%%%%%%%%%%%%%%%%%%%%%%%%%%%%%%%%%%%%%%%%%%%%%%%%%%%%%%%%%%%%%%%%%%%%%%%%%%%%
%%%%%%%%%%%%%%%%%%%%%%%%%%%%%%%%%%%%%%%%%%%%%%%%%%%%%%%%%%%%%%%%%%%%%%%%%%%%%%%%%%%%

\section{Appendices}
\label{sec:appendix}

\setcounter{theorem}{0}
\setcounter{equation}{0}

%\begin{section}\label{sec:appendix}

Here we will present additional calculations, which are mainly based on the papers~\cite{A84},~\cite{A86},~\cite{stream},~\cite{persson-lic},
clarifying the theory being used to prove Lemma~\ref{derivatalemma}.
We begin with deriving the $p$-Laplace equation~\eqref{sepeq} in polar coordinates, which brings us to the \emph{separation equation}~\eqref{sepeq2}. Then we will develop a \emph{stream function} technique in order to handle the situation when $1<p<2$.

%\begin{appendix}

\subsection{Transforming the $p$-laplacian to polar coordinates}\label{transform}

%\subsubsection*{A1: Transforming to polar coordinates}

The $p$-Laplace equation \eqref{eq:plapequation} can be transformed to polar coordinates by putting
$x(r,\phi)=r \cos \phi$, $y(r,\phi)=r\sin\phi$
and hence $u(x,y)=u[x(r,\phi),y(r,\phi)]$.
%In this paper we study growth of a $p$-harmonic measure, defined below, which is a generalization of harmonic measure, related to the $p$-Laplace equation. For $p\in(1,\infty)$, the $p$-Laplace equation yields
%
%\begin{align*}\label{eq:plapequation}
%\Delta_{p} u :=\nabla \cdot \left( |\nabla u |^{p - 2} \nabla  u \right) = 0,
%\end{align*}
%in divergence form. We will now transform this equation to polar coordinates.
%undergraduate vector calculus gives
Introduce  $\psi=|\nabla u|^2$ and note that when $\psi \neq 0$ the equation is equivalent to
%
%\begin{multline*}
%\Delta_{p} u :=\nabla \cdot \left( |\nabla u |^{p - 2} \nabla  u \right)=\nabla \cdot \left(|\nabla u|^{p-2}\right)\cdot\nabla u+
%    |\nabla u|^{p-2}\nabla ^2 u =\\
%    (p-2)\left(|\nabla u| ^{p-3}\nabla(|\nabla u|)\right)\cdot\nabla u+|\nabla u|^{(p-2)}\nabla^2u=\\
%    |\nabla u|^{(p-2)}\left(\nabla^2u+\frac{(p-2)}{|\nabla u|}\nabla(|\nabla u|)\cdot \nabla u\right)=0.
%\end{multline*}
%
%Thus
%\begin{equation}
%\label{lap}
%\nabla^2u+\frac{(p-2)}{|\nabla u|}\nabla\left(|\nabla %u|\right)\cdot \nabla u=0
%\end{equation}
%Also
%\begin{equation*}
%\nabla (|\nabla u|^2)=2|\nabla u|\nabla(|\nabla u|)\Rightarrow \nabla(|\nabla u|)=\frac{1}{2|\nabla u|}\nabla(|u|^2)
%\end{equation*}
%inserted in equation~\ref{lap} to obtain
%
%\begin{equation*}
%  \nabla^2u+\frac{(p-2)}{2|\nabla u|^2}\nabla\left(|\nabla u|^2\right)\cdot\nabla u=0.
%\end{equation*}
%
%Define a help function $\psi$ by $\psi=|\nabla u|^2$, which,  gives
%
\begin{equation*}
   \nabla^2u+\frac{(p-2)}{2\psi}\nabla\psi\cdot\nabla u=0.
\end{equation*}
%
%To simplify complicated things lets transform this form of the $p$-Laplace equation to a more convenient version in polar coordinates i.e.
%
%Put $x(r,\phi)=r \cos \phi$, $y(r,\phi)=r\sin\phi$ and hence $u(x,y)=u[x(r,\phi),y(r,\phi)]$.
Trivial calculations yield $u_r=\partial_r u=\frac{\partial u}{\partial r}=\cos \phi \frac{\partial u}{\partial x}+\sin \phi \frac{\partial u}{\partial y}$ and %similarly
$u_\phi=\partial_\phi u=\frac{\partial u}{\partial \phi}=-r\sin \phi \frac{\partial u}{\partial x}+r\cos \phi \frac{\partial u}{\partial y}$. %First year linear algebra simplifies otherwise long and tedious calculations.
Put
\begin{equation*}
P =
\begin{pmatrix}
\cos\phi & -\sin\phi  \\
\sin\phi & \cos\phi  \\
\end{pmatrix}
%\end{equation*}
\hspace{1cm}\text{giving}\hspace{1cm}
P^{-1}=P^{T}=
\begin{pmatrix}
\cos\phi & \sin\phi  \\
-\sin\phi & \cos\phi  \\
\end{pmatrix},
\end{equation*}
%
%since $P^{T}P=I$,
and thus in operator matrix notation
\begin{equation*}
\begin{pmatrix}
\partial_r \\
\frac{1}{r}\partial_\phi \\
\end{pmatrix}
=
%\begin{pmatrix}
%\cos\phi & \sin\phi  \\
%-\sin\phi & \cos\phi  \\
%\end{pmatrix}
P^T
\begin{pmatrix}
\partial_x \\
\partial_y \\
\end{pmatrix}
\hspace{1cm}\text{and}\hspace{1cm}
\begin{pmatrix}
\partial_x \\
\partial_y \\
\end{pmatrix}
=
%\begin{pmatrix}
%\cos\phi & -\sin\phi \\
%\sin\phi &\cos\phi\\
%\end{pmatrix}
P
\begin{pmatrix}
\partial_r \\
\frac{1}{r}\partial\phi
\end{pmatrix}.
\end{equation*}
%
%The helpfunction $\psi$ can be calculated by
It follows that
\begin{equation*}
\psi=|\nabla u|^2
=\left(P\begin{pmatrix}
\partial_r u \\
\frac{1}{r} \partial_{\phi}u\\
\end{pmatrix}
\right)^T
P
\begin{pmatrix}
\partial_r u\\
\frac{1}{r}\partial_{\phi}
\end{pmatrix}
=P^TP\left(\partial_r u,\frac{1}{r}\partial_{\phi}u\right)
\begin{pmatrix}
\partial_r u\\
\frac{1}{r}\partial_{\phi}
\end{pmatrix}
=u_{r}^2+\frac{1}{r^2}u_{\phi}^2,
\end{equation*}
%
%Also
giving $\psi_r=2u_ru_{rr}-\frac{2}{r^3}u_\phi^2+\frac{2}{r^2}u_\phi u_{r\phi}$ and $\psi_\phi=2u_ru_{r\phi}+\frac{2}{r^2}u_\phi u_{\phi\phi}.$
Therefore
\begin{align*}
 \nabla\psi\cdot\nabla u=\left(P\begin{pmatrix}
\partial_r \psi \\
\frac{1}{r} \partial_{\phi}\psi\\
\end{pmatrix}
\right)^T
P
\begin{pmatrix}
\partial_r u\\
\frac{1}{r}\partial_{\phi}
\end{pmatrix}
=P^TP\left(\partial_r \psi,\frac{1}{r}\partial_{\phi}\psi\right)
\begin{pmatrix}
\partial_r u\\
\frac{1}{r}\partial_{\phi}
\end{pmatrix}
=\psi_ru_r+\frac{1}{r^2}\psi_{\phi}u_{\phi}\\
=2u_r^2u_{rr}-\frac{2}{r^3}u_ru_\phi^2+\frac{2}{r^2}u_ru_\phi u_{r\phi}+
\frac{2}{r^2}u_ru_{\phi}u_{r\phi}+\frac{2}{r^4}u_{\phi}^2u_{\phi \phi}.
\end{align*}
%
%The helpfunction $\psi$ also needs be differentiated with respect of $r$ and $\phi$ respectively, hence
% \begin{equation*}
% \psi_r=\frac{\partial \psi}{\partial r}=2u_ru_{rr}-\frac{2}{r^3}u_\phi^2+\frac{2}{r^2}u_\phi u_{r\phi}
% \hspace{0.5cm}
%\text{and}\hspace{0.5cm}
% \psi_\phi=\frac{\partial \psi}{\partial \phi}=2u_ru_{r\phi}+\frac{2}{r^2}u_\phi u_{\phi \phi}.
% \end{equation*}
%Therefore
%
%\begin{equation*}
%\nabla\psi\cdot\nabla u=%\psi_ru_r+\frac{1}{r^2}\psi_{\phi}u_{\phi}=
%2u_r^2u_{rr}-\frac{2}{r^3}u_ru_\phi^2+\frac{2}{r^2}u_ru_\phi u_{r\phi}+
%\frac{2}{r^2}u_ru_{\phi}u_{r\phi}+\frac{2}{r^4}u_{\phi}^2u_{\phi \phi}.
%\end{equation*}
%
Recalling the Laplace operator in polar coordinates,
%Straight forward but long and tedious calculations gives the ``ordinary'' $(p=2)$ Laplace differential operator in polar   coordinates
%
%\begin{multline*}
$$
\Delta_{(r,\phi)}=%\frac{\partial^2}{\partial x^2}+\frac{\partial^2}{\partial y^2}=\left(\cos\phi \frac{\partial}{\partial_r}-\frac{\sin\phi}{r}\frac{\partial}{\partial\phi}\right)
%\left(\cos\frac{\partial}{\partial_r}-\frac{\sin\phi}{r}\frac{\partial}{\partial\phi}\right)+\\ \left(\sin\phi \frac{\partial}{\partial_r}+\frac{\cos\phi}{r}\frac{\partial}{\partial\phi}\right)
%\left(\sin\frac{\partial}{\partial_r}+\frac{\cos\phi}{r}\frac{\partial}{\partial\phi}\right)=
\frac{\partial^2}{\partial r^2}+\frac{1}{r}\frac{\partial}{\partial r}+\frac{1}{r^2}\frac{\partial^2}{\partial\phi^2},
$$
%\end{multline*}
%
we obtain, for $1  < p  < \infty$,
%Assume that $(p<1\leq\infty)$ and $p\neq 2$, hence

\begin{multline*}
\nabla^2u+\frac{(p-2)}{2\psi}\nabla\psi\cdot\nabla u=u_{rr}+\frac{1}{r}u_r+\frac{1}{r^2}u_{\phi\phi}\\  +\frac{(p-2)}{2(u_r^2+\frac{1}{r}u_{\phi}^2)}\left(2u_r^2u_{rr}-\frac{2}{r^3}u_ru_{\phi}^2+
\frac{2}{r^2}u_ru_\phi u_{r\phi}+\frac{2}{r^2}u_ru_\phi u_{r\phi}+\frac{2}{r^4}u_\phi^2u_{\phi\phi}\right)=0.
\end{multline*}

\bigskip

\noindent
Put, when $p\neq 2$, $p-2=1/b$, multiply by $2\psi=2(u_r^2+\frac{1}{r^2}u_{\phi}^2)$ and simplify to finally arrive at the $p$-Laplace equation in polar coordinates:
\begin{multline}\label{sepeq2}
\left(b+1\right)u_r^2u_{rr}+\frac{b}{r^2}\left(u_{rr}u_{\phi}^2+u_r^2u_{\phi\phi}\right)+ \frac{(b+1)}{r^4}u_{\phi}^2u_{\phi\phi}+\frac{b}{r}u_r^3+\frac{(b-1)}{r^3}u_ru_{\phi}^2+\frac{2}{r^2}u_ru_\phi u_{r\phi}=0.\tag{\ref{eq:polar-coords}}
\end{multline}

%\end{appendix}

%\begin{appendix}
\subsection{Searching for solutions of the form $u(r,\phi)=r^kf(\phi)$}\label{search}

%\subsubsection*{A2: Searching for solutions of the form $u(r,\phi)=r^kf(\phi)$}

We are searching for solutions of the form $u(r,\phi)=r^kf(\phi)$, where  $f(\phi)\in \mathcal{C}^2$ and $k$ are to be determined. %Also we have to assume that $f(\phi)\in \mathcal{C}^2$ and that $r>0$ of course.
Differentiation with respect to $r$ and $\phi$ yields $u_r=kr^{k-1}f(\phi$), $u_{rr}=k(k-1)r^{k-2}f(\phi)$, $u_\phi=r^kf'(\phi)$ and $u_{\phi\phi}=r^kf''(\phi)$, which inserted in Equation~\eqref{eq:sepeqinsec3} yield %and simplifying gives

\begin{equation}\label{sepeq}
\left[\left(b+1\right)(f')^2+bk^2f^2\right]f''+\left(2k+bk-1\right)kf(f')^2+\left(bk+k-1\right)k^3f^3=0,\tag{\ref{eq:sepeqinsec3}}
\end{equation}

\bigskip

\noindent
which is our \emph{separation equation}~\eqref{eq:sepeqinsec3}. We are initially interested in the case where $p>2$, $f(\phi)>0$ and $k>0$, which is case $\alpha$ in~\cite{A86}.
Following~\cite{A86} there are three cases to consider $ak>1$, $ak=1$ and $0<ak<1$. Here $a=(p-1)/(p-2)$. First, we have a look at the case $ak>1$ and later (if needed) the cases $ak=1$, $0<ak<1$. The equation~(\ref{sepeq2}) can be solved (for details check~\cite[pages 141-143]{A86}) to yield

\bigskip

\begin{equation}\label{losning}
\left\{
	\begin{array}{ll}
		\phi=\theta^{\ast}+\displaystyle \int_{\theta^{\ast}}^\theta\frac{a-\cos^2\theta'}{ak-\cos^2\theta'}\,d\theta'  \\[24pt]
		f(\phi)=c\left(1-\frac{\cos^2\theta(\phi)}{ak}\right)^\frac{k-1}{2}\cos\theta(\phi)
	\end{array}
\right.
\end{equation}

\bigskip

\noindent
for some constant $c>0$. For convenience further ahead we put $c=(\frac{ak-1}{ak})^{-\frac{k-1}{2}}$ giving $f_{\nu,p}(0)=1$. Note that $\phi(\theta)$ is monotonously increasing since $\phi'(\theta)=\frac{a-\cos^2\theta}{ak-\cos^2\theta}>0$, which holds since $ak>1$ and $a > 1$. %the denominator can never be less then zero since $ak>1$ and the nominator can never be less than zero because $a=(p-1)/(p-2)>1$, for $2<p<\infty$.
Hence there exists an inverse $\theta=\phi^{-1}$, which is continuous since $\phi(\theta)$ is continuous.

\bigskip

%\end{appendix}

\subsection{Finding the radial exponent $k=k(\nu,p)$}\label{find}

To find the radial exponent $k=k(\nu,p)$ we need to solve the parametric equation~\eqref{losning} and from that determine $k$.
To do so put $\theta^{\ast}=0$ in~\eqref{losning} and observe that
\begin{equation*}
\frac{a-\cos^2\theta'}{ak-\cos^2\theta'}=\frac{ak-\cos^2\theta'-a(k-1)}{ak-\cos^2\theta'}=1-\frac{a(k-1)}{ak-\cos^2\theta'}\,,
\end{equation*}
%
%\bigskip
%
which transforms the integrand to
\begin{equation*}
\phi=\theta-a(k-1)\int_0^{\theta}\frac{d\theta'}{ak-\cos^2\theta'}\,.
\end{equation*}
Using partial fraction decomposition yields
%Partial fraction decomposition is appropriate to use here, thus
%
\begin{equation*}
\phi=\theta-a(k-1)\int_0^{\theta}\frac{d\theta'}{ak-\cos^2\theta'}=\theta-\frac{a(k-1)}{2\sqrt{ak}}\left(\,\int_0^{\theta}\frac{d\theta'}{\sqrt{ak}+\cos\theta'}+\int_0^{\theta}\frac{d\theta'}{\sqrt{ak}\cos\theta'}\,\right).
\end{equation*}
%
%\bigskip
%
Utilizing the $\tan(\theta/2)$ substitution and simplifying we arrive at
\begin{multline*}%\label{arctan}
\phi=\theta-\frac{a(k-1)}{\sqrt{ak}\sqrt{ak-1}}\left[\arctan\left(\lambda\tan\frac{\theta}{2}\right)+\arctan\left(\frac{1}{\lambda}\tan\frac{\theta}{2}\right)\right]\\
=\theta-(1-\frac{1}{k})\frac{\sqrt{ak}}{\sqrt{ak-1}
}\left[\arctan\left(\lambda\tan\frac{\theta}{2}\right)+\arctan\left(\frac{1}{\lambda}\tan\frac{\theta}{2}\right)\right]\,,
\end{multline*}
where
\begin{equation*}
\lambda=\frac{\sqrt{ak-1}}{\sqrt{ak}+1}\hspace{1cm}\text{and}\hspace{1cm}\frac{1}{\lambda}=\frac{\sqrt{ak-1}}{\sqrt{ak}-1}.
\end{equation*}
%
%\bigskip
%
Define
\begin{multline*}
\bar{\phi}=\phi\left(\frac{\pi}{2}\right)-\phi\left(-\frac{\pi}{2}\right)=\frac{\pi}{2}-(1-\frac{1}{k})\frac{\sqrt{ak}}{\sqrt{ak-1}
}\left[\arctan\left(\lambda\tan\frac{\pi}{4}\right)+\arctan\left(\frac{1}{\lambda}\tan\frac{\pi}{4}\right)\right]\\
-\left(\frac{-\pi}{2}-(1-\frac{1}{k})\frac{\sqrt{ak}}{\sqrt{ak-1}
}\left[\arctan\left(\lambda\tan(-\frac{\pi}{4})\right)+\arctan\left(\frac{1}{\lambda}\tan(-\frac{\pi}{4})\right)\right]\right)\\
=\pi-(1-\frac{1}{k})\frac{2\sqrt{ak}}{\sqrt{ak-1}
}\left[\arctan\left(\lambda\right)+\arctan\left(\frac{1}{\lambda}\right)\right].
\end{multline*}
%
%\bigskip
%
Now $\frac{d}{d\lambda}\left(\arctan(\lambda)+\arctan(1/\lambda)\right)=\frac{1}{1+\lambda^2}-\frac{1}{1+\lambda^2}=0$, and hence $\arctan(\lambda)+\arctan(1/\lambda)$ is constant.
Therefore $\lambda=\pm 1$ determines the function values for $\lambda>0$ and $\lambda<0$ respectively, so that

\begin{equation*}
\arctan\left(\lambda\right)+\arctan\left(\frac{1}{\lambda}\right)=\left\{
	\begin{array}{ll}
		\frac{\pi}{2} &,\quad \lambda \geq0 \\[10pt]
		-\frac{\pi}{2} &,\quad\lambda <0.\\
	\end{array}
\right.
\end{equation*}
Since $ak>1$ we must have $\lambda>0$ and therefore $\bar{\phi}$ becomes
\begin{equation*}
  \bar{\phi}=\phi\left(\frac{\pi}{2}\right)-\phi\left(-\frac{\pi}{2}\right)=\pi\left(1-\left(1-\frac{1}{k}\right)\frac{\sqrt{ak}}{\sqrt{ak-1}}\right).
\end{equation*}
%
%\bigskip
%
Note that $\phi\left(\frac{\pi}{2}\right)=\frac{\bar{\phi}}{2}$,  $\phi\left(-\frac{\pi}{2}\right)=-\frac{\bar{\phi}}{2}$, and also
$f\left(\frac{\bar\phi}{2}\right)=f\left(-\frac{\bar\phi}{2}\right)=0$ which suits our purpose perfectly.

Recall the sector $\mathcal{S}_\nu$ from \eqref{eq:def-S_v}, having aperture $\pi/\nu$ and apex at the origin.
In $\mathcal{S}_\nu$ we can now introduce our continuous $p$-harmonic function $u_{\nu,p}(x,y)=r^kf_{\nu,p}(\phi)$, where $f_{\nu,p}$ can be written as
\begin{equation}\label{losningv}
f_{\nu,p}(\phi)=c\left(1-\frac{\cos^2\theta_{\nu,p}(\phi)}{ak}\right)^{\frac{k-1}{2}}\cos\theta_{\nu,p}(\phi),\tag{\ref{eq:losningv}}
\end{equation}

\bigskip

\noindent
where, for convenience, we take $c=(\frac{ak-1}{ak})^{-\frac{k-1}{2}}$. Note that since $ak>1$, $f_{\nu,p}(\phi)$ is bounded by a constant, depending only on $\nu$ and $p$, % $c(v,p)$
when $\phi\in(-\frac{\pi}{2\nu},\frac{\pi}{2\nu})$.

\bigskip

When $|\phi|<\frac{\pi}{2\nu}$, we have that
\begin{equation*}
  \phi=\theta_{\nu,p}(\phi)-\left(1-\frac{1}{k}\right)\frac{\sqrt{ak}}{\sqrt{ak-1}}\left[\arctan\left(\lambda_{\nu,p}\tan\frac{\theta_{\nu,p}(\phi)}{2}\right)+\arctan\left(\frac{1}{\lambda_{\nu,p}}\tan\frac{\theta_{\nu,p}(\phi)}{2}\right)\right].
\end{equation*}
%
%\bigskip
%
The condition that determines the radial exponent $k = k(\nu,p)$ is given by
\begin{equation*}
  \frac{\pi}{\nu}=\phi\left(\frac{\pi}{2}\right)-\phi\left(-\frac{\pi}{2}\right)=\pi\left(1-\left(1-\frac{1}{k}\right)\frac{\sqrt{ak}}{\sqrt{ak-1}}\right).
\end{equation*}
%
%\bigskip
%
Recalling that $a=(p-1)/(p-2)$ and solving for $k$ we obtain two roots
\begin{equation}\label{radialexponent}
k_1(\nu,p)=\frac{\sqrt{(1-2\nu)(p-2)^2+\nu^2p^2}(\nu-1)+(2-p)(1-2\nu)+\nu^2p}{2(p-1)(2\nu-1)}
\end{equation}
and
\begin{equation*}
k_2(\nu,p)=\frac{\sqrt{(1-2\nu)(p-2)^2+\nu^2p^2}(1-\nu)+(2-p)(2\nu-1)+\nu^2p}{2(p-1)(2\nu-1)}.
\end{equation*}

\bigskip

\noindent
To decide which of these two solutions to choose for $k$ we put $p=2$ giving $k_1(\nu,2)=\nu$ and $k_2(\nu,2)=\nu/(\nu-1)$. Therefore, $k(\nu,p)=k_1(\nu,p)$ is the true solution (it matches $k(\nu,2) = \nu$ and $k_2(\nu,2)$ fails to be positive).

Differentiating~\eqref{radialexponent} with respect to $\nu$ gives
%\begin{equation*}
%    \lim_{p\rightarrow\infty}ak_1(v,p)=\frac{v^2+2v-1+(v-1)\sqrt{(v-1)^2}}{2(2v-1)},
%\end{equation*}
%similarly
%\begin{equation*}
%    \lim_{p\rightarrow\infty}ak_2(v,p)=\frac{(v-1)\left((v-1)-\sqrt{(v-1)^2}\right)}{2(2v-1)}.
%\end{equation*}

%Since $v^2>v-1/2$ for all $v\in\mathbb{R}$ we only have to consider the case $ak>1$.

\begin{equation*}%\label{radialv}
     \frac{\partial k(\nu,p)}{\partial \nu}= \nu\, \frac{p(\nu-1)\sqrt{(\nu-1)^2p^2+4(2\nu-1)(p-1)}+(\nu-1)^2p^2+2(2\nu-1)(p-1)}{(p-1)(2\nu-1)^2\sqrt{(\nu-1)^2p^2+4(2\nu-1)(p-1)}},
\end{equation*}

\bigskip
\noindent
which is %Equation~\eqref{radialv} is
easily seen to be greater than zero if $\nu\geq 1$. For $\nu\in[\frac{1}{2},1)$ $\frac{\partial k}{\partial \nu}$ is nonnegative if $(\nu-1)^2p^2\geq -p(\nu-1)\sqrt{(\nu-1)^2p^2+(8\nu-4)p-8\nu+4}$, which leads us to the inequality $4(p-1)(2\nu-1)\geq0$. Hence $\frac{\partial k}{\partial \nu}\geq0$, for all $p\in(1,\infty)$ and $\nu\in[\frac{1}{2},\infty)$.
When $p=\infty$ the conclusion follows by differentiation on \eqref{eq:kinfhej}.

 Given the expression~\eqref{radialexponent} for $k$ we can now check if the case $ak > 1$ gives the desired $p$-harmonic function whenever $\nu \in [\frac{1}{2},\infty)$ and $p \in (2,\infty)$. Since $k$ is increasing in $\nu$ we only have to check the worst case scenario i.e. $ak(\frac{1}{2},p)=\frac{p-1}{p-2}\frac{p-1}{p}=1+\frac{1}{p^2-2p}>1$, for all $p>2$.
 Therefore, it is enough to consider the case $ak > 1$.

Differentiating~\eqref{radialexponent} with respect to $p$ gives

\begin{equation}\label{radialp}
        \frac{\partial k(\nu,p)}{\partial p}=(1-\nu)\,\frac{(\nu-1)\sqrt{(\nu-1)^2p^2+4(2\nu-1)(p-1)}+\nu^2p+(2\nu-1)(p-2)}{2(2\nu-1)(p-1)^2\sqrt{(\nu-1)^2p^2+4(2\nu-1)(p-1)}}.
\end{equation}

\bigskip

\noindent
 Putting $\frac{\partial k}{\partial p}=0$ gives either $\nu=1$ or $(\nu-1)\sqrt{(\nu-1)^2p^2+4(2\nu-1)(p-1)}+\nu^2p+(2\nu-1)(p-2)=0$. The latter holds when $\nu=0$ (which is not allowed), $\nu=1/2$ or when $p=1$ (of which none are allowed). Going to the limit in~\eqref{radialp} yields $\lim_{\nu\to\frac{1}{2}}\frac{\partial k}{\partial p}=\frac{1}{p^2}>0$, for all $p>1$.
 %(note that since $k(\frac{1}{2},p)=\frac{p-1}{p}$ is differentiable for all \komN{Also strange as it is indep. of $v$:} $v\geq 1/2$, therefore $\partial k(1/2,p)/\partial p=1/p^2$).
 Since $\frac{\partial k}{\partial p}$ is zero only when $\nu=1$ it is sufficient to investigate two points $\nu_1\in[1/2,1)$ and $\nu_2\in(1,\infty)$ in order to know the sign of the derivative. If $\nu_1=3/4$ then the numerator of~\eqref{radialp} is equal to $\frac{\sqrt{p^2+32p-32}+17p-16}{64}>0$ for all $p > 1$. Similarly, if $\nu_2=2$ then the numerator of~\eqref{radialp} is equal to $6-7p-\sqrt{p^2+12p-12}<0$ for all $p\geq 4\sqrt{3}-6$ and thus also for $p>1$.
 We conclude that $k(\nu,p)$ is increasing in $p$ for all $\nu\in[1/2,1)$, decreasing in $p$ for $\nu>1$ and constant if $\nu=1$ ($k(1,p)=1$).

  \bigskip

\subsection{Computing the derivative $f'(\phi)$}\label{compute}

We also need an estimation of the derivative of $f(\phi)$. Differentiation of $f$ in Equation~\eqref{losning} and  simplifying yields

\begin{align}\label{uppskattad derivata}
\frac{df}{d\theta}=c\left(1-\frac{\cos^2\theta}{ak}\right)^{\frac{k-3}{2}}\left(\frac{\cos^2\theta}{a}-1\right)\sin\theta.
\end{align}

\bigskip

\noindent
Since $\theta^{\ast}=0$ we have $\phi(\theta)= \int_{0}^\theta\frac{a-\cos^2\theta'}{ak-\cos^2\theta'}\,d\theta'$, and differentiation gives
%$\theta=\phi^{-1}\left( \int_{0}^\theta\frac{a-\cos^2\theta'}{ak-\cos^2\theta'}\,d\theta'\right)$ and %this inverse exist due to the monotonicity. Implicit differentiation gives

\begin{align}\label{implicit derivata}
\frac{d\phi}{d\theta}=\frac{d}{d\theta}\left( \int_{0}^\theta\frac{a-\cos^2\theta'}{ak-\cos^2\theta'}\,d\theta'\right)=\frac{a-\cos^2\theta}{ak-\cos^2\theta}.
\end{align}

\bigskip

\noindent
 Now $\frac{df}{d\theta}=\frac{d f}{d\theta}\frac{d\theta}{d\phi}$ by the chain rule and recall the fact that $\phi$ is monotone, continuous, differentiable and hence invertible. For simplicity let $\phi=g(\theta)$ be the right hand side of the integrand in~\eqref{losning} so that $\theta=g^{-1}(\phi)$. By the inverse function theorem \cite[Theorem~9.24]{rudin}  $\frac{dg^{-1}}{d\phi}$ exists and once again by the chain rule $1=\frac{dg^{-1}}{d\phi}\frac{d\phi}{d\theta}=\frac{d\theta}{d\phi}\frac{d\phi}{d\theta}$.
 Thus $\frac{d\theta}{d\phi}=(\frac{d\phi}{d\theta})^{-1}=\frac{ak-\cos^2\theta}{a-\cos^2\theta}$ and by using ~\eqref{uppskattad derivata} with~\eqref{implicit derivata} we arrive at

\begin{align}\label{bounded derivative}
    \frac{df}{d\phi}=\frac{d f}{d\theta}\frac{d\theta}{d\phi}=c\left(1-\frac{\cos^2\theta}{ak}\right)^{\frac{k-3}{2}}\left(\frac{\cos^2\theta}{a}-1\right)\frac{ak-\cos^2\theta}{a-\cos^2\theta}\sin\theta
    =-ck\left(1-\frac{\cos^2\theta}{ak}\right)^{\frac{k-1}{2}}\sin\theta.\tag{\ref{eq:bounded derivative}}
\end{align}

\subsection{A little taste of stream functions}\label{taste}

Here we will present a simple stream function technique for partial differential equations of the form

\begin{equation*}
    \nabla\cdot\left(\frac{F(|\nabla u|)}{|\nabla u|}\nabla u\right)=0,
\end{equation*}

\bigskip

\noindent
where $F(t)>0$ is monotonically increasing and continuously differentiable on $(\alpha,\beta)\in\Omega\subset\mathbb{R}^2$. Applying this technique to the $p$-harmonic equation will reveal a $q$-harmonic \emph{stream function}, where $p$ and $q$ are conjugate exponents ($1/p+1/q=1$). This has been described earlier in~\cite{stream} but for the convenience of the reader we will give a short presentation of it here.

\bigskip

Consider $\Omega\subset\mathbb{R}^2$, $u\in C^2(\Omega)$ and assume also that $0\leq\alpha<|\nabla u |<\beta$ in $\Omega$ for constants $\alpha$ and $\beta$. Now if

\begin{equation}\label{stream1}
    \nabla\cdot\left(\frac{F(|\nabla u|)}{|\nabla u|}\nabla u\right)=0,
\end{equation}
%
%i.e.
it follows that

\begin{equation*}
    \frac{\partial}{\partial x}\left(\frac{F(|\nabla u|)}{|\nabla u|}\frac{\partial u}{\partial x}\right)=    \frac{\partial}{\partial y}\left(-\frac{F(|\nabla u|)}{|\nabla u|}\frac{\partial u}{\partial y}\right).
\end{equation*}

\bigskip

\noindent
Define $\psi_y=\left(\frac{F(|\nabla u|)}{|\nabla u|}u_x\right)$ and $\psi_x=\left(-\frac{F(|\nabla u|)}{|\nabla u|}u_y\right)$.  Since $F(t),u_x\in C^1(\Omega)$ and $|\nabla u|>0$ it follows that $\psi_x$ and $\psi_y$ are integrable, thus $\psi=\int\psi_y\,dy+C(x)$ and $\psi=\int\psi_x\,dx+D(y)$, for some functions $C(x)$ and $D(y)$. Further since $|\nabla u|>0$ both $\frac{\partial}{\partial x}(|\nabla u|)$ and $\frac{\partial}{\partial y}(|\nabla u|)$ exists, hence $\psi\in C^2(\Omega)$.

Now $|\nabla\psi|=\frac{|F(|\nabla u|)|}{|\nabla u|}(u_x^2+u_y^2)^{\frac{1}{2}}=F(|\nabla u|)$ and $|\nabla u|=F^{-1}(|\nabla\psi|)$, since $F(t)$ is strictly increasing and hence invertible, therefore
\begin{equation*}
    u_x=\frac{F^{-1}(|\nabla\psi|)}{|\nabla \psi|}\psi_y\quad\text{and}\quad u_y=-\frac{F^{-1}(|\nabla\psi|)}{|\nabla \psi|}\psi_x.
\end{equation*}

\noindent
Now, since $u_{xy}=u_{yx}$ we deduce
\begin{multline}\label{stream2}
\nabla\cdot\left(\frac{F^{-1}(|\nabla\psi|)}{|\nabla \psi|}\nabla\psi\right)=\frac{F^{-1}(|\nabla\psi|)}{|\nabla \psi|}\left(\frac{\partial}{\partial x}\hat{x}+\frac{\partial}{\partial y}\hat{y}\right)\left(\frac{\partial\psi}{\partial x}\hat{x}+\frac{\partial\psi}{\partial y}\hat{y}\right)\\
=\frac{F^{-1}(|\nabla\psi|)}{|\nabla \psi|}\left(\frac{\partial^2\psi}{\partial x \partial y}-\frac{\partial^2\psi}{\partial y \partial x}\right)=0
\end{multline}

\bigskip

\noindent
in $\Omega$. Conversely if we begin with $\psi$ satisfying~\eqref{stream2} and proceed similarly we will find $u$ satisfying~\eqref{stream1}. Equations~\eqref{stream1} and~\eqref{stream2} is said to constitute a \emph{reciprocal pair} of equations. Also

\begin{multline*}
    \nabla u\cdot\nabla \psi= \frac{F^{-1}(|\nabla\psi|)}{|\nabla \psi|}\left(\frac{\partial\psi}{\partial y}\hat{x}-\frac{\partial\psi}{\partial x}\hat{y}\right)\cdot\left(\frac{\partial\psi}{\partial x}\hat{x}+\frac{\partial\psi}{\partial y}\hat{y}\right)\\
    =
    \frac{F^{-1}(|\nabla\psi|)}{|\nabla \psi|}\left(\frac{\partial\psi}{\partial y}\frac{\partial\psi}{\partial x}-\frac{\partial\psi}{\partial x}\frac{\partial\psi}{\partial y}\right)=0,
\end{multline*}

\bigskip

\noindent
so the gradient of $\psi$ is perpendicular to the streamlines of $u$. Thus streamlines of $u$ are level curves of $\psi$ and vice versa. In fluid mechanics $\psi$ is called the \emph{stream function} corresponding to the \emph{potential} $u$ (or conversely), see~\cite{bers}. Let $1<p<\infty$ and consider $\nabla \cdot ( |\nabla u |^{p - 2} \nabla  u ) = 0$ so that $F(t)=t^{p-1}$. We then have $t=F^{-1}(s)=s^{\frac{1}{p-1}}$, and the corresponding reciprocal equation
\begin{equation*}
   \nabla\cdot\left(|\nabla\psi|^{\frac{2-p}{p-1}}\, \nabla\psi\right)=0.
\end{equation*}

\bigskip

\noindent
Since $p$ and $q$ are conjugate exponents we have $q=p/(p-1)$ and the reciprocal equation becomes
\begin{equation*}
   \nabla\cdot\left(|\nabla\psi|^{q-2}\, \nabla\psi\right)=0.
\end{equation*}

\bigskip
\noindent
Thus, the reciprocal of the $p$-harmonic equation in the plane is the $q$-harmonic equation, where $1/p+1/q=1$. The above discussion will lead to the result below (it is sometimes presented as a definition), which is proven in~\cite{stream}:

\begin{lemma}\label{stream}
Let $p\in(1,\infty)$ and let $u$ be $p$-harmonic ($u$ not constant) in a simply connected domain $\Omega\subset\mathbb{R}^2$. Then there exists a $q$-harmonic function $v\in C^1(\Omega)$, where $1/p+1/q=1$, such that

\begin{equation*}
\left\{
	\begin{array}{ll}
		v_x=-|\nabla u|^{p-2}u_y \\[24pt]
		v_y=|\nabla u|^{p-2}u_x.
	\end{array}
\right.
\end{equation*}

\bigskip
\noindent
Both $u$ and $v$ have locally Hölder continuous gradients. The zeroes of $\nabla u$ and $\nabla v$ are isolated in $\Omega$. Streamlines of $u$ are level lines of $v$ and vice versa.
\end{lemma}
For any vector, define operator $T:\mathbb{R}^2\rightarrow\mathbb{R}^2$ by\;
$T\left(\begin{array}{ccc}
a \\
b \\
\end{array}\right)=
\left(\begin{array}{ccc}
-b \\
a \\
\end{array}\right)$. According to Lemma~\ref{stream} we then have

\begin{equation}\label{cooloperator}
\nabla v=|\nabla u|^{p-2}\,T(\nabla u)=|\nabla u|^{p-1}\,T\left(\frac{\nabla u}{|\nabla u|}\right).
\end{equation}

 \bigskip
 \noindent
 Hence $|\nabla v|=|\nabla u|^{p-1}$ and $|\nabla v|^q=|\nabla u|^p$. Let us find stream functions $v$ to our radial function $u=r^kf(\phi)$ in polar coordinates. The below result is established in~\cite{stream} and a complex version can be found in~\cite{persson-lic}.

 \begin{lemma}\label{konjugatlemma}
 Let $u(r,\phi)=r^kf(\phi)$ be $p$-harmonic in the sector $\mathcal{S}_\nu$, $k>0$, and suppose that $p\in(2,\infty)$. Then there exists a $q$-harmonic stream function $v(\lambda,\phi)=r^{\lambda}g(\phi)$, where $\lambda=(p-1)(k-1)+1$, $q=p/(p-1)$, and

 \begin{equation*}%\label{streamphi}
     g(\phi)=-\dfrac{1}{\lambda}f'(\phi)\left(k^2f(\phi)^2+f'(\phi)^2\right)^{\frac{p-2}{2}}.%\tag{\ref{eq:streamphi}}
 \end{equation*}

 \bigskip

\noindent
The function $g(\phi)$ is periodic whenever $f(\phi)$ is.
\end{lemma}

\noindent
\begin{proof}
Assume $u(r,\phi)=r^k f(\phi)$ to be $p$-harmonic in $\mathcal{S}_\nu$, and $k>0$. The nabla operator in polar coordinates gives
\begin{equation*}
    \nabla u=\begin{pmatrix}
\cos\phi & -\sin\phi  \\
\sin\phi & \cos\phi  \\
\end{pmatrix}
\begin{pmatrix}
\partial_r \\
\frac{1}{r}\partial_\phi \\
\end{pmatrix}
r^k f(\phi)=r^{k-1}\left(k f(\phi)e(\phi)+f'(\phi)d(\phi)\right),
\end{equation*}

\bigskip

\noindent
where $e(\phi)=\begin{pmatrix}
\cos\phi   \\
\sin\phi   \\
\end{pmatrix}$
and $d(\phi)=\begin{pmatrix}
-\sin\phi   \\
\cos\phi   \\
\end{pmatrix}$. Now, substituting $|\nabla u|=r^{k-1}\sqrt{\left(k^2f(\phi)^2+f'(\phi)^2\right)}$ into~\eqref{cooloperator}, we obtain

\begin{multline*}
    |\nabla u|^{p-1}T\left(\frac{\nabla u}{|\nabla u|}\right)=r^{(k-1)(p-1)}\left(k^2f(\phi)^2+f'(\phi)^2\right)^{\frac{p-1}{2}}T\left(\frac{kf(\phi)e(\phi)+f'(\phi)d(\phi)}{\sqrt{k^2f(\phi)^2+f'(\phi)^2}}\right)\\
    %=\frac{r^{(k-1)(p-1)}\left(k^2f(\phi)^2+f'(\phi)^2\right)^{\frac{p-1}{2}}}{\sqrt{k^2f(\phi)^2+f'(\phi)^2}}T\begin{pmatrix}
%kf(\phi)\cos\phi-f'(\phi)\sin\phi   \\
%kf(\phi)\sin\phi+f'(\phi)\cos\phi   \\
%\end{pmatrix}\\
%=\frac{r^{(k-1)(p-1)}\left(k^2f(\phi)^2+f'(\phi)^2\right)^{\frac{p-1}{2}}}{\sqrt{k^2f(\phi)^2+f'(\phi)^2}}
%\begin{pmatrix}
%-kf(\phi)\sin\phi-f'(\phi)\cos\phi   \\
%kf(\phi)\cos\phi-f'(\phi)\sin\phi   \\
%\end{pmatrix}\\
=\frac{r^{(k-1)(p-1)}\left(k^2f(\phi)^2+f'(\phi)^2\right)^{\frac{p-1}{2}}}{\sqrt{k^2f(\phi)^2+f'(\phi)^2}}\left(-f'(\phi)e(\phi)+kf(\phi)d(\phi)\right).
\end{multline*}

\bigskip

\noindent
Now we search for a stream function of the form
$v=r^{\lambda}g(\phi)$ such that
$$
\nabla v=r^{\lambda-1}\left(\lambda g(\phi)e(\phi)+g'(\phi)d(\phi)\right)
$$
holds.
It is convenient to separate the direction from the modulus, i.e.,

\begin{equation*}
    \nabla v=r^{\lambda-1}\sqrt{\lambda^2\left(g(\phi)^2+g'(\phi)^2\right)}\;\frac{\lambda g(\phi)e(\phi)+g'(\phi)d(\phi)}{\sqrt{\lambda^2\left(g(\phi)^2+g'(\phi)^2\right)}}.
\end{equation*}

\bigskip

\noindent
Equation~\eqref{cooloperator} and the fact that $|\nabla v|=|\nabla u|^{p-1}$ give the following system of equations

\begin{equation*}
\left\{
	\begin{array}{ll}
		r^{\lambda-1}\sqrt{\lambda^2g(\phi)^2+g'(\phi)^2}=r^{(k-1)(p-1)}\left(k^2f(\phi)^2+f'(\phi)^2\right)^{\frac{p-1}{2}}, \\[18pt]
		\frac{\lambda g(\phi)}{\sqrt{\lambda^2g(\phi)^2+g'(\phi)^2}}=-\frac{f'(\phi)}{\sqrt{k^2f(\phi)^2+f'(\phi)^2}},\\[18pt]
		\frac{g'(\phi)}{\sqrt{\lambda^2 g(\phi)^2+g'(\phi)^2}}=\frac{k f(\phi)}{\sqrt{k^2f(\phi)^2+f'(\phi)^2}}.
	\end{array}
\right.
\end{equation*}

\bigskip
\noindent
Hence $\lambda-1=(p-1)(k-1)$ and the other conditions yield

\begin{equation}\label{sepeq3}
\left\{
	\begin{array}{ll}
		\lambda^2g(\phi)^2+g'(\phi)^2=\left(k^2f(\phi)^2+f'(\phi)^2\right)^{p-1}, \\[18pt]
		\lambda g(\phi)=-f'(\phi)\left(k^2f(\phi)^2+f'(\phi)^2\right)^{\frac{p-2}{2}},\\[18pt]
		g'(\phi)=k f(\phi)\left(k^2f(\phi)^2+f'(\phi)^2\right)^{\frac{p-2}{2}}.
	\end{array}
\right.
\end{equation}

\bigskip
\noindent
Clearly $\lambda=(p-1)(k-1)+1$ and $g(\phi)=-\dfrac{1}{\lambda}f'(\phi)\left(k^2f(\phi)^2+f'(\phi)^2\right)^{\frac{p-2}{2}}$.
We remark that a straightforward calculation verifies that the system \eqref{sepeq3} is identical to the separation equation~\eqref{eq:sepeqinsec3}, which $f(\phi)$ is known to satisfy.
\end{proof}

%\bigskip
%The system of equations in \eqref{sepeq3} gives an ordinary differential equation for $f(\phi)$ which is  in fact identical to the separation equation~\eqref{eq:sepeqinsec3}, which $f(\phi)$ is known to satisfy. To verify this is straightforward and left as an exercise for the reader.

\bigskip

%%%%%%%%%%%%%%%%%%%%%%%%%%%%%%%%%%%%%%%%%%%%%%%%%%%%%%%%%%%%%%%%%%%%%%%%%%%%%%%%%%%%%%
%%%%%%%%%%%%%%%%%%%%%%%%%%%%%%%%%%%%%%%%%%%%%%%%%%%%%%%%%%%%%%%%%%%%%%%%%%%%%%%%%%%%%%
%%%%%%%%%%%%%%%%%%%%%%%%%%%%%%%%%%%%%%%%%%%%%%%%%%%%%%%%%%%%%%%%%%%%%%%%%%%%%%%%%%%%%%
%%%%%%%%%%%%%%%%%%%%%%%%%%%%%%%%%%%%%%%%%%%%%%%%%%%%%%%%%%%%%%%%%%%%%%%%%%%%%%%%%%%%%%
%%%%%%%%%%%%%%%%%%%%%%%%%%%%%%%%%%%%%%%%%%%%%%%%%%%%%%%%%%%%%%%%%%%%%%%%%%%%%%%%%%%%%%
%%%%%%%%%%%%%%%%%%%%%%%%%%%%%%%%%%%%%%%%%%%%%%%%%%%%%%%%%%%%%%%%%%%%%%%%%%%%%%%%%%%%%%

\noindent
{\bf Acknowledgement.}
We would like to thank Marcus Olofsson and an anonymous reviewer for valuable comments and suggestions on an earlier version of this  manuscript.
This work was partially supported by the Swedish research council grant 2018-03743.

\end{document}